\documentclass[11pt]{amsart}

\usepackage{latexsym}
\usepackage{amsmath,amsthm,amssymb, amscd,color}
\usepackage{mathrsfs}
\usepackage[all]{xy}
\usepackage[all]{xy}
\usepackage{tikz}
\usepackage{tikz-cd}
\usepackage{mathtools}
\usetikzlibrary{snakes}
\usepackage{blkarray} 
\usepackage{framed}

\usepackage{anysize}
\marginsize{3cm}{3cm}{2.5cm}{2.5cm}

\theoremstyle{plain}
\newtheorem{theorem}{Theorem}[section]
\newtheorem{lemma}[theorem]{Lemma}
\newtheorem{corollary}[theorem]{Corollary}
\numberwithin{equation}{section}

\theoremstyle{definition}
\newtheorem{definition}[theorem]{Definition}
\newtheorem{example}[theorem]{Example}

\newtheorem{proposition}[theorem]{Proposition}
\newtheorem{remark}[theorem]{Remark}

\theoremstyle{remark}

\newcommand{\CC}{\mathbb{C}}

\newcommand{\RR}{\mathbb{R}}
\newcommand{\ZZ}{\mathbb{Z}}

\newcommand{\eqcont}[1]{%
	\mathrel{\ooalign{\hbox{\scalebox{1.5}[1]{$\rhd$}}\cr%
	\kern-0.1ex\raise1.15ex\hbox{\scalebox{1.5}[1]{$\sim$}}\cr%
	\kern1.5ex\raise-1.2ex\hbox{\scalebox{0.7}{#1}}\cr%
	}}}

\begin{document}

\title[Various topological complexities]{Various topological complexities of small covers and real Bott manifolds}

\author[K. Brahma]{Koushik Brahma}
\address{Department of Mathematics, Indian Institute of Technology Madras, Chennai 600036, India}
\email{koushikbrahma95@gmail.com}

\author[B. Naskar]{Bikramaditya Naskar}
\address{Department of Mathematics, Indian Institute of Technology Madras, Chennai 600036, India}
\email{bikramadityaix@gmail.com}

\author[S. Sarkar]{Soumen Sarkar}
\address{Department of Mathematics, Indian Institute of Technology Madras, Chennai 600036, India}
\email{soumen@iitm.ac.in}

\author[S. Sau]{Subhankar Sau}
\address{Indian Statistical Institute, Kolkata-700108, India}
\email{subhankarsau18@gmail.com}

\date{\today}
\thanks{
}

\abstract In this paper, we compute the LS-category and equivariant LS-category of a small cover and its real moment angle manifold. We calculate a tight lower bound for the topological complexity of many small covers over a product of simplices. Then we compute symmetric topological complexity of several small covers over a product of simplices. We calculate the LS one-category of real Bott manifolds and infinitely many small covers.
\endabstract

\maketitle 

\section{Introduction}

The topological complexity of a space is a numerical homotopy invariant introduced by Michael Farber in \cite{Far03}, which connects motion planning problems in robotics. Briefly, given a mechanical system $\mathcal{M}$, a motion planning algorithm for $\mathcal{M}$ is a function that associates to any pair of states $(a, b)$ of $\mathcal{M}$ to a continuous motion of the system starting at $a$ and ending at $b$. Interestingly, the topological complexity is a particular case of another homotopy invariant called the `sectional category' of a map $p \colon E \rightarrow B$ where $E$ and $B$ are path connected spaces. The sectional category of $p$, denoted by $\mbox{secat}(p)$, is the least integer $k$ such that there is an open cover $\lbrace U_1,..., U_k \rbrace$ of $B$, and there is a local section $s_i \colon U_i \rightarrow E$ of $p$ for each $i$ satisfying $p \circ s_i= id_{U_i} \colon U_i \hookrightarrow B$ where $id_{U_i}$ denotes the inclusion. We remark that the genus of a fibration was introduced by Schwarz \cite{Sva58}. However, James \cite{Jam78} used `sectional category' instead of `genus'.

Let $Y$ be the space of all possible configurations of a mechanical system. We assume that $Y$ is a Hausdorff path-connected topological space. Let $PY$ be the space of all continuous paths $\gamma \colon [0,1] \rightarrow Y$ in $Y$ equipped with the compact-open topology. 
 Consider the path fibration 
 \begin{equation} \label{tc_first_equation}
 \pi \colon PY \rightarrow Y \times Y
 \end{equation}
defined by $\pi(\gamma)=( \gamma(0), \gamma(1))$. A motion planning algorithm of $Y$ is defined by a section $s \colon Y \times Y \rightarrow PY$ of the fibration $\pi$. This section exists if and only if $Y$ is contractible. Interestingly in general, almost all configuration spaces are non-contractible. To compute the complexity of motion-planning algorithm for a non-contractible space $Y$, Farber defined the topological complexity of $Y$ by the sectional category of $\pi$. The survey \cite{Far-Lec08} contains several introductory results related to motion planning.

A symmetric version of the topological complexity arises when one restricts the local planners for which the motion from $a$ to $b$ is the reverse of the motion from $b$ to $a$ and the motion from $a$ to $a$ is constant. In notation, consider a map $s \colon Y \times Y \rightarrow PY$ (not necessarily continuous) such that  $\pi \circ s= \textit{Id}_{Y \times Y}$ and $s(a, a)(t)= a$, ${s(a, b)(t)= s(b, a)(1-t)}$ for all $a, b \in Y$ and $t \in [0, 1]$. This motivates the notion of symmetric topological complexity, given by Farber and Grant in \cite{FG07}. Some developments in symmetric topological complexity can be found in \cite{GoPe, Gon, Gra}.

Consider a continuous partial section $s \colon U \rightarrow PY$  of the fibration $\pi$ over an open subset $U \subseteq Y \times Y$. The map $s$ can be described as a homotopy $h \colon U \times [0, 1] \rightarrow Y$ defined by $h(u, t)= s(u)(t)$ for $u \in U, t \in [0, 1]$. Let $p_1 \colon Y \times Y \rightarrow Y$ and $p_2 \colon Y \times Y \rightarrow Y$ denote the projections onto the first and the second factor, respectively. Since $s$ is a section, the homotopy $h$ connects $h(u, 0)= p_1(u)$ and $h(u, 1)= p_2(u)$.
Therefore the open sets $U_i \subseteq Y \times Y$, which appear in the definition of topological complexity, can be equivalently characterized by the property that their two projections $U_i \rightarrow Y$ on the first and the second factors are homotopic. For an aspherical space $Y$, a connected subspace $U$ of $Y$ which is homotopy equivalent to a cell complex, the set of homotopy classes of maps $U \rightarrow Y$ is in a one-to-one correspondence with the set of conjugacy classes of homomorphisms $\pi_1(U, u_0) \rightarrow \pi_1(Y, y_0)$. Using this idea, Farber, Grant, Lupton, and Oprea introduced $\mathbf{TC}^{\mathcal{D}}(Y)$, the $\mathcal{D}$-topological complexity for a path-connected topological space, see  \cite{FGLO19bredon}. Here the letter `$\mathcal{D}$' in the notation $\mathbf{TC}^{\mathcal{D}}(Y)$ stands for the `diagonal'. In \cite{FGLO19}, Farber, Grant, Lupton, and Oprea introduced some properties of $\mathcal{D}$-topological complexity. Note that symmetric topological complexity is not homotopy invariant but $\mathcal{D}$-topological complexity is homotopy invariant. Some related results can be found in \cite{Dra}.

 A small cover of dimension $n$ is an $n$-dimensional closed smooth manifold with a locally standard $\ZZ_2^n$-action whose orbit space is a simple polytope. It was introduced in the pioneering paper \cite{DJ91} as a generalization of real projective toric varieties. An $n$-dimensional toric variety is an algebraic normal variety that admits an action of $(\CC^*)^n$ with an open dense orbit. A non-singular complete toric variety is simply called a toric manifold. The real locus of a toric manifold is called a real toric manifold. A real Bott tower is a sequence of smooth complete real toric varieties, see Subsection \ref{subsec_gen_real_bott_cohomo_ring}. In this paper, we compute lower and upper bounds for the topological complexity, symmetric topological complexity, and LS one-category of a class of small covers and real Bott manifolds.

The paper is organized as follows. In Section \ref{sec_cohomo_ring_small_cover_real_bott}, we study the definition of small cover over a simple polytope, generalized real Bott manifold, and the relation between them. We modify the cohomology ring of a small cover over a product of simplices $\prod_{j=1}^m\Delta^{n_j}$ with $\ZZ_2$ coefficients as $\ZZ_2[y_1, y_2, \dots, y_m]/I$ where $I$ is given in \eqref{compactform_cohomo_ring_ideal}. We prove $y_j^{n_j} \neq 0$ in the cohomology ring $H^*(M^n(P, \lambda); \ZZ_2)$ as in \eqref{compactform_cohomo_ring}, see Lemma \ref{yj_nj_neq_0}. We also recall the notion of real moment angle manifolds and complexes. 

In Section \ref{sec_eqlscat_small_cover}, we recall the definition and some properties of LS-category and equivariant LS-category of a topological space. We compute the LS-category and equivariant LS-category of a small cover. We calculate the LS-category of the real moment angle manifold for $r$-gon and the equivariant LS-category of a real moment angle complex.

In Section \ref{sec_higher_tc_small_cover}, we give a tight lower bounds to the topological complexity of a small cover over a product of two simplices. We compute the topological complexity for some classes of real Bott manifolds. 

In Section \ref{sec_symm_tc_small_cover}, we rewrite the definition and some basic properties of symmetric topological complexity and give bounds for the symmetric topological complexity of several small covers over a product of simplices. 

Finally in Section \ref{sec_dtc_small_cover}, we recall the definition and some basic properties of the LS one-category and $\mathcal{D}$-topological complexity. We calculate the exact value of LS one-category of a simple polytope when its real moment angle manifold is simply connected and orientable. We calculate LS one-category of a small cover over a product of simplices, and give bounds of $\mathcal{D}$-topological complexity for a small cover over a product of simplices.



\section{Cohomology rings of small covers, generalized real Bott manifolds, and real moment angle complexes} \label{sec_cohomo_ring_small_cover_real_bott}


In this section, we recall simple polytopes and the constructive definition of a small cover over a simple polytope using \cite{DJ91}. If the polytope is a product of finitely many simplices, then the small cover is known as a generalized real Bott manifold. We give a presentation of the cohomology ring of a generalized real Bott manifold. Later, we study real moment angle manifolds and complexes.

\subsection{Small covers and its cohomology ring}\label{subsec_small_cover_cohomo_ring}

In this subsection, we recall the definition of small cover and its cohomology ring with $\ZZ_2$-coefficients following \cite{DJ91}.

A convex polytope is a convex hull of finitely many points in $\RR^n$ for some $n \in \ZZ_{\geq 0}$. 
The face of dimension $0$ and $(n-1)$ in a convex polytope of dimension $n$ are called the vertex and the facet of the polytope, respectively. The vertex set and the facet set of a convex polytope $P$ are denoted by $V(P)$ and $\mathcal{F}(P)$, respectively.
An $n$-dimensional convex polytope is called simple if at each vertex exactly $n$ many facets intersect. Throughout this paper, we denote an $n$-dimensional simple polytope by $P$.

\begin{definition} \label{Def 2.1}
A function $\lambda \colon \mathcal{F}(P) \rightarrow \mathbb{Z}_2^n $ is called a characteristic function if the submodule of $\mathbb{Z}_2^n$ generated by $\lbrace \lambda(F_{i_1}),..., \lambda(F_{i_\ell}) \rbrace$ is an $\ell$-dimensional direct summand of $\mathbb{Z}_2^n$ whenever $F_{i_1} \cap \dots \cap F_{i_\ell} \neq \emptyset$.
The vector $\lambda_i:= \lambda(F_i)$ is called the characteristic vector associated with the facet $F_i$ for $i=1,..., r$, and the pair $(P, \lambda)$ is called a characteristic pair.
\end{definition}

We recall the construction of a small cover from a characteristic pair $(P, \lambda)$.
For each point $p \in P$, let $F(p)$ be the unique face of $P$, which contains $p$ in its relative interior. Let $F(p)= F_{i_1} \cap \cdots \cap F_{i_k}$ for some unique facets $F_{i_1},..., F_{i_k}$. Define $G_{F(p)}$ as a subgroup of $\mathbb{Z}_2^n$ generated by $\lambda(F_{i_1}),...,\lambda(F_{i_k})$. We define an equivalence relation on $P \times \mathbb{Z}_2^n$ as follows:
$$(p,g) \backsim (q,h) \Leftrightarrow p=q, g^{-1}h \in G_{F(p)}.$$

The identification space $M^n(P,\lambda): =( P\times \ZZ_2^n ) / \sim$ has an $n$-dimensional manifold structure with a natural $\ZZ_2^n$-action induced by the group operation on the second factor of $P\times \ZZ_2^n$. The projection onto the first factor gives the orbit map
\begin{equation}\label{orbit map}
\rho \colon M^n(P,\lambda)  \to P \text{ defined by } [p,g]_{\sim} \mapsto p \nonumber,
\end{equation}
where $[p,g]_{\sim}$ is the equivalence class of $(p,g)$. The manifold $M^n(P, \lambda)$ is called a small cover over $P$ with the characteristic function $\lambda$, see \cite{DJ91} for details.

  Let $\lbrace F_1,..., F_r \rbrace$ be the facets of $P$ and the indeterminates $v_1,...,v_r$ correspond  bijectively to the facets $F_1,...,F_r$ respectively.
\begin{proposition}\label{cohomoring_smallcover_polyt} \cite [Theorem 4.14]{DJ91}
Let $\rho \colon M^n(P, \lambda) \rightarrow P$ be a small cover over a simple polytope $P$ with $|\mathcal{F}(P)|=r$. Then $$H^*(M^n(P, \lambda), \mathbb{Z}_2) \cong \mathbb{Z}_2[v_1,...,v_r]/(\Tilde{I}+\Tilde{J}),$$
 where the ideal $\Tilde{I}$ is generated by the monomials $v_{s_1} \cdots v_{s_\ell}$, if $F_{s_1} \cap \cdots \cap F_{s_\ell}= \emptyset,$ and the ideal $\Tilde{J}$ is generated by the $n$ coordinates of the vector $\Lambda_{\Tilde{J}}$ where $\Lambda_{\Tilde{J}}= \sum_{i=1}^{r} \lambda_iv_i$.
 \end{proposition}


\begin{example}
The $n$-dimensional real projective space $\mathbb{RP}^n$ is an example of a small cover over the $n$-dimensional simplex $\Delta^n$. A finite product of $\mathbb{RP}^n$'s is also a small cover.
\end{example}



\subsection{Generalized real Bott manifolds and its cohomology ring} \label{subsec_gen_real_bott_cohomo_ring}


In this subsection, we study generalized real Bott manifolds and give a nice presentation of its cohomology ring with $\ZZ_2$-coefficients.

A generalized real Bott tower of height $m$ is a sequence 
\begin{equation} \label{bott_tower}
   B_m \xrightarrow{\pi_m} B_{m-1} \xrightarrow{\pi_{m-1}} \cdots \xrightarrow{\pi_2} B_1 \xrightarrow{\pi_1} B_0= \lbrace \mbox{pt} \rbrace  
\end{equation}
of manifolds $B_j= \mathbb{P}(\underline{\mathbb{R}} \oplus E_j^{(1)} \oplus \cdots \oplus E_j^{(n_j)} )$, where $\underline{\mathbb{R}}$ is the trivial line bundle over $B_{j-1}$, $E_j^{(i)}$ is a real line bundle over $B_{j-1}$ for $i=1,..., n_j$, and $j= 1,..., m$. Here $\mathbb{P}(.)$ denotes the projectivization. The space $B_j$ is called a $j$-th stage generalized real Bott manifold. In this case, when $n_j= 1$ for every $j$, $B_j$ is called a real Bott manifold.

\begin{proposition}  \cite[Corollary 4.6]{KZ16}
The $j$-th stage generalized real Bott manifold $B_j$ of the tower \eqref{bott_tower} is a small cover over $\prod_{i=1}^j \Delta^{n_i}$ where $\Delta^{n_i}$ is the $n_i$-simplex. 
\end{proposition}
The converse statement also holds by the following proposition.
\begin{proposition} \cite[Proposition 2.7]{DU19} \label{small_cover_realbott} 
    Every small cover over a product of simplices is a generalized real Bott manifold.
\end{proposition}

Now we discuss the cohomology ring of a small cover over a finite product of simplices. Let 
\begin{equation}\label{P is a product of simplices}
P:= \prod_{j=1}^m \Delta^{n_j},
\end{equation}
where $\Delta^{n_j}$ is a simplex of dimension $n_j$. Then, the dimension of $P$ is $n := \sum_{j=1}^m n_j$. 
Let
\begin{equation}\label{Eq_define Ns}
\mathcal{N}_s:= \sum_{j=1}^s n_j,
\end{equation}
for $s=1, \dots, m$. Thus $\mathcal{N}_1=n_1$ and $\mathcal{N}_m=n$. Let us assume $\mathcal{N}_0:=0$.

Let $V(\Delta^{n_j}):=\{v_0^j, \dots ,v_{n_j}^j\}$ be the vertices of $\Delta^{n_j}$ for $j=1, \dots,m$. Then the vertex set of $P$ is given by
\begin{equation}\label{Eq_vertex set of P}
V(P):=\{v_{\ell_1 \ell_2 \dots \ell_m}:=(v_{\ell_1}^1, v_{\ell_2}^2, \dots, v_{\ell_m}^m )~ |~0\leq \ell_j \leq n_j\}.
\end{equation} 
Let $\mathcal{F}(\Delta^{n_j}):=\{F_0^{\Delta_j}, \dots, F_{n_j}^{\Delta_j}\}$ be the facets of $\Delta^{n_j}$ where the facet $F_{k_j}^{\Delta_j}$ does not contain the vertex $v_{k_j}^j$ for $j=1, \dots, m$. So, the facet set of $P$ is
\begin{equation}\label{Eq_facet set of P}
\mathcal{F}(P):=\{F_{k_j}^j ~|~0\leq k_j \leq n_j, j=1, \dots, m\},
\end{equation}
where $F_{k_j}^j:=\Delta^{n_1} \times \dots \times \Delta^{n_{j-1}} \times F_{k_j}^{\Delta_j} \times \Delta^{n_{j+1}}\times \dots \times \Delta^{n_m}.$
Observe that 
the vertex $v_{\ell_1 \ell_2 \dots \ell_m}$ is the unique intersection of the $n$-many facets of $\mathcal{F}(P) \setminus \{F_{\ell_j}^j ~|~j=1, \dots,m\}$. In particular,
\begin{equation}\label{at_v00000}
v_{0 \dots 0}= F_1^1 \cap \dots \cap F_{n_1}^1 \cap \dots \cap F_1^m \cap \dots \cap F_{n_m}^m.
\end{equation}

Let
\begin{equation}\label{Eq_define_lambda_on_P}
\lambda \colon \mathcal{F}(P) \to \mathbb{Z}_2^n
\end{equation}
be a $\mathbb{Z}_2$-characteristic function on $P$ where $P$ is the product of simplices as in \eqref{P is a product of simplices}. Then from \eqref{at_v00000}, we have $\{\lambda(F^1_1), \dots, \lambda(F^1_{n_1}), \dots,\lambda(F^m_1), \dots, \lambda(F^m_{n_m})\}$ is a basis of $\ZZ_2^n$ over $\mathbb{Z}_2$.
So, we may assume that these vectors are assigned with the standard basis vectors. 
Thus,
 \begin{align*} 
\lambda(F_1^j)=\boldsymbol{e_{\mathcal{N}_{j-1}+1}}, & \dots, \lambda(F_{n_j}^j)=\boldsymbol{e_{\mathcal{N}_j}},
\end{align*}
for $j=1,..., m$.
The remaining $m$ facets $\lbrace F_0^1,..., F_0^m \rbrace$ are assigned with the vectors as follows
 \begin{equation} \label{Eq_lambda_to_other_facets}
\lambda(F_0^j):=\boldsymbol{\alpha_j} \in \mathbb{Z}_2^n \quad ~\text{for}~j=1, \dots, m,
 \end{equation}
so that the above assignment satisfies Definition \ref{Def 2.1}. This gives us vector matrices of order $(1 \times m)$ and $(m \times m)$, and a scalar matrix of order $(n \times m)$ as following:
 \begin{align}\label{Eq_vector matrix}
  A:=&
\begin{pmatrix} 
  \boldsymbol{\alpha_1} & \boldsymbol{\alpha_2} & \dots & \boldsymbol{\alpha_m}
\end{pmatrix}_{1 \times m}
=
\begin{pmatrix}
 \boldsymbol{\alpha_1^1} & \dots & \boldsymbol{\alpha_m^1} \\
\vdots & \dots & \vdots \\
\boldsymbol{\alpha_1^m} & \dots & \boldsymbol{\alpha_m^m} 
\end{pmatrix}_{m \times m}=
\begin{pmatrix}
\alpha_{11}^1 & \dots & \alpha_{m1}^1\\
\vdots & \dots & \vdots \\
\alpha^1_{1n_1} & \dots & \alpha^1_{mn_1}\\
\vdots & \dots & \vdots \\
\alpha^m_{11} & \dots & \alpha^m_{m1} \\
\vdots & \dots & \vdots \\
 \alpha^m_{1n_m} & \dots & \alpha^m_{mn_m} 
\end{pmatrix}_{n \times m}, \nonumber
\end{align} 
where $\boldsymbol{\alpha_j} \in \mathbb{Z}_2^n$ is the $j$-th column vector of $A$, $\boldsymbol{\alpha_j^k} \in \mathbb{Z}_2^{n_k}$ is the $(k,j)$-th entry of the $m \times m$ vector matrix and $\alpha_{ji}^k \in \ZZ_2$ is the $(\mathcal{N}_{k-1}+i,j)$-th entry of the $n \times m$ scalar matrix.
Throughout this paper, the vectors $\boldsymbol{e_i}$ and $\boldsymbol{\alpha_j}$ of $\ZZ_2^n$ are considered as the column entries of the matrices for $i=1, \dots, n$, and $j=1, \dots, m$.
 

Now we calculate the cohomology ring of the small cover $M^n(P, \lambda)$ when $P$ is a product of simplices as in \eqref{P is a product of simplices} and the characteristic function $\lambda$ on $P$ is given by \eqref{Eq_define_lambda_on_P}.
Let us assign the indeterminate $x_i$ to the facet $F^j_{k_j}$ where
 $$i=(\sum_{s=1}^{j-1} n_s )+ k_j=\mathcal{N}_{j-1}+k_j,$$ 
for $1 \leq k_j \leq n_j$, $j=1, \dots, m$. Therefore $i \in \lbrace 1, \dots, n \rbrace$.
We also assign the indeterminate $x_i$ to the facet $F^j_0$ where $i=n+j$ for $j=1, \dots, m$. Note that $F^j_1 \cap \dots \cap F^j_{n_j} \cap F^j_0 = \emptyset$. Then, from Proposition \ref{cohomoring_smallcover_polyt}, we have 
\begin{equation}\label{Eq_cohom over small cover}
H^*(M^n(P,\lambda); \ZZ_2) \cong \ZZ_2[x_1, \dots, x_{n+m}]/ \Tilde{I} + \Tilde{J},
\end{equation}
where the ideals $\Tilde{I}$ and $\Tilde{J}$ are as follows. The ideal $\Tilde{I}$ is given by 
\begin{equation}\label{Eq_the ideal I}
\Tilde{I}= \big< \{ x_{\mathcal{N}_{j-1}+1} x_{\mathcal{N}_{j-1}+2} \dots x_{\mathcal{N}_j} x_{n+j} ~|~ j=1, \dots, m \}\big>,
\end{equation} 
where $\mathcal{N}_j$ is defined in \eqref{Eq_define Ns}.
The ideal $\Tilde{J}$ is generated by the coordinates of
\small{\begin{equation}\label{eq_lamda_j}
   \Lambda_{\Tilde{J}}=\begin{pmatrix}
   \lambda(F_1)^t & \lambda(F_2)^t & \dots & \lambda(F_{n+m})^t
   \end{pmatrix}_{(n \times (n+m))} \cdot \begin{pmatrix}
   x_1 & x_2 & \dots & x_{n+m}
   \end{pmatrix}^t_{(n+m) \times 1}.
\end{equation}}
\noindent In \eqref{eq_lamda_j}; for $i=1, \dots, n$, we denote $F_i=F^j_{k_j}$ with $i=\mathcal{N}_{j-1}+k_j$ for $1 \leq k_j \leq n_j$, $j=1, \dots, m$ and for $i=n+1, \dots, n+m$, we denote $F_i=F^j_0$ with $i=n+j$ where $j=1, \dots, m$.
  
  Note that $\Lambda_{\Tilde{J}}$ is an $n$ tuple. The $i$-th coordinate of $\Lambda_{\Tilde{J}}$ is
  $$x_i+\alpha_{1k_j}^j x_{n+1} + \alpha_{2k_j}^j x_{n+2} + \dots + \alpha^j_{mk_j} x_{n+m},$$
  where $i=\mathcal{N}_{j-1}+k_j$; $k_j=1, \dots, n_j$ and $j=1, \dots, m$.
  Thus any $x_i$ can be written as a $\ZZ_2$-linear combination of $x_{n+1}, \dots, x_{n+m}$ for $i=1, \dots, n$. For simplicity, we denote the indeterminate $x_{n+j}$ by $y_j$ for $j=1, \dots, m$. Thus,
\begin{equation} \label{x_i_relation_with_y_i}
  x_i= \sum_{\ell=1}^m \alpha^j_{\ell k_j} y_{\ell} ~\mbox{where}~ i=\mathcal{N}_{j-1}+k_j, k_j=1, \dots, n_j ~\mbox{and}~ j=1, \dots, m,
\end{equation} 
  in $H^*(M^n(P, \lambda); \ZZ_2)$. Then the generators of the ideal $\Tilde{I}$ in \eqref{Eq_the ideal I} can be described in terms of $y_j$'s. Therefore, we have 
  \begin{equation} \label{compactform_cohomo_ring}
      H^*(M^n(P, \lambda); \ZZ_2) \cong \ZZ_2[y_1, y_2, \dots, y_m]/I, \text{where}
  \end{equation}
\begin{equation} \label{compactform_cohomo_ring_ideal}
I=\big< \{ \prod_{k_j=1}^{n_j} \big( \sum_{\ell=1}^m \alpha_{\ell k_j}^j y_{\ell} \big) y_j ~|~ j=1, \dots, m \} \big>.
\end{equation}

We have the following observation on the cohomology ring.

\begin{lemma}\label{yj_nj_neq_0}
   Let  $M^n(P, \lambda)$ be a small cover over a finite product of simplices with the characteristic function $\lambda$ as in \eqref{Eq_define_lambda_on_P}. Then, for $j \in \lbrace 1,...,m \rbrace$, $y_j^{n_j} \neq 0$ in the cohomology ring $H^*(M^n(P, \lambda); \ZZ_2)$ as in \eqref{compactform_cohomo_ring}.
\end{lemma}

\begin{proof}

Let $P:= \prod_{j=1}^m \Delta^{n_j}$ be a product of $m$ simplices. We know the cohomology ring of a small cover over a product of simplices from \eqref{compactform_cohomo_ring}.
The function $\lambda$ determines the following $m \times m$ vector matrix $A$.

\begin{align}
  A:=&
\begin{pmatrix}
\boldsymbol{\alpha_1^1} & \boldsymbol{\alpha_2^1} & \dots & \boldsymbol{\alpha_m^1} \\
\boldsymbol{\alpha_1^2} & \boldsymbol{\alpha_2^2} & \dots & \boldsymbol{\alpha_m^2} \\
\vdots & \vdots & \dots & \vdots \\
\boldsymbol{\alpha_1^m} & \boldsymbol{\alpha_2^m} & \dots & \boldsymbol{\alpha_m^m} \\
\end{pmatrix}_{m \times m}. \nonumber
\end{align}

Therefore, by the arguments in \cite[Proposition $5.1$]{CMS10}, $A$ is conjugate to a unipotent lower triangular vector matrix of the following form:

\begin{align}\label{lower_triangular_matrix}
\Tilde{A}:=&
\begin{pmatrix}
\textbf{1} & \boldsymbol{0} & \dots & \boldsymbol{0} \\
\boldsymbol{\beta_1^2} & \textbf{1} & \dots & \boldsymbol{0} \\
\vdots & \vdots & \dots & \vdots \\
\boldsymbol{\beta_1^m} & \boldsymbol{\beta_2^m} & \dots & \textbf{1} \\
\end{pmatrix}_{m \times m}, 
\end{align}
where $\boldsymbol{\beta_j^k}= ( \beta_{j1}^k, \beta_{j2}^k, \dots, \beta_{jn_k}^k )^t \in \mathbb{Z}_2^{n_k}$ and $\textbf{1}=(1, \dots, 1)^t \in \ZZ^{n_k}_2$ for $k=1, \dots, m$. The matrix $\Tilde{A}$ is called the Bott matrix. Thus the ideal $\Tilde{J}$ is generated by the coordinates of the following matrix.
\begin{equation*}
\left(\begin{array}{cccccccccccc}
1 & 0 & \cdots & 0 & 1 & 0 & \cdots & 0 & \cdots & \cdots & 0 & 0\\
0 & 1 & \cdots & 0 & 1 & 0 & \cdots & 0 & \cdots & \cdots & 0 & 0 \\
\vdots & \vdots & \ddots & \vdots & \vdots & \vdots & \ddots & \vdots & \ddots & \ddots & \vdots & \vdots \\
0 & 0 & \cdots & 1 & 1 & 0 & \cdots & 0 & \cdots & \cdots & 0 & 0 \\
0 & 0 & \cdots & 0 & \beta_{11}^2 & 1 & \cdots & 1 & \cdots & \cdots & 0 & 0 \\
\vdots & \vdots & \ddots & \vdots & \vdots & \vdots & \ddots & \vdots & \ddots & \ddots & \vdots & \vdots \\
0 & 0 & \cdots & 0 & \beta_{1n_2}^2 & 0 & \cdots & 1 & \cdots & \cdots & 0 & 0 \\
\vdots & \vdots & \ddots & \vdots & \vdots & \vdots & \ddots & \vdots & \ddots & \ddots & \vdots & \vdots \\
\vdots & \vdots & \ddots & \vdots & \vdots & \vdots & \ddots & \vdots & \ddots & \ddots & \vdots & \vdots \\
0 & 0 & \cdots & 0 & \beta_{11}^m & 0 & \cdots & \beta_{21}^m & \cdots & \cdots & 0 & 1 \\
\vdots & \vdots & \ddots & \vdots & \vdots & \vdots & \ddots & \vdots & \ddots & \ddots & \vdots & \vdots \\
0 & 0 & \cdots & 0 & \beta_{1n_m}^m & 0 & \cdots & \beta_{2n_m}^m & \cdots & \cdots & 1 & 1 \\
\end{array}\right)_{n \times (n+ m)}
\begin{pmatrix}
x_1 \\
x_2 \\
\vdots\\
x_{\mathcal{N}_1} \\
y_1\\
x_{\mathcal{N}_1 +1} \\
\vdots \\
y_2\\
\vdots \\
\vdots \\
\vdots\\
x_{\mathcal{N}_m}\\
y_m\\
\end{pmatrix}_{(n+m) \times 1}.
\end{equation*}

 Let $\alpha_j:= x_{\mathcal{N}_{j-1}+1} x_{\mathcal{N}_{j-1}+2} \cdots x_{\mathcal{N}_j} y_j $ for $j=1,...,m$. Here $\alpha_j$'s are generators of the ideal $I$ in \eqref{compactform_cohomo_ring}. From the above matrix multiplication, the first $n_1$ elements are $x_1+ y_1= 0$, $x_2+ y_1= 0$,..., $x_{\mathcal{N}_1}+ y_1= 0$. Therefore we get $x_1= x_2= \cdots= x_{\mathcal{N}_1}= y_1$.\\
So,
$$\alpha_1= x_1 x_2 \cdots x_{n_1} y_1= y_1^{n_1+1}.$$



We have the following using \eqref{x_i_relation_with_y_i}.
 \begin{align*}
\alpha_j &= x_{\mathcal{N}_{j-1}+1} x_{\mathcal{N}_{j-1}+2} \cdots x_{\mathcal{N}_j} y_j \\
&= (y_j+ \beta_{11}^jy_1+ \beta_{21}^jy_2+ \cdots + \beta_{(j-1)1}^jy_{j-1})(y_j+ \beta_{12}^jy_1+ \beta_{22}^jy_2+ \\
& \cdots + \beta_{(j-1)2}^jy_{j-1}) 
 \cdots (y_j+ \beta_{1n_j}^jy_1+ \beta_{2n_j}^jy_2+ \cdots + \beta_{(j-1)n_j}^jy_{j-1})y_j,
\end{align*}
for $j=2,...,m$. Now the least power of $y_j$ in $\alpha_j$ is $n_j+1$. Our claim is that $y_j^{n_j} \neq 0$. If not, let $y_j^{n_j}= 0$. Then $y_j^{n_j} \in I$. But the least power of $y_j$ which appears as a term in a polynomial in the ideal $I$ is $y_j^{n_j+1}$. This is a contradiction. Hence $y_j^{n_j} \notin I$, i.e., $y_j^{n_j} \neq 0$ in $H^*(M^n(P, \lambda); \mathbb{Z}_2)$ for $j=1,2,..., m$.

\end{proof}

\subsection{Real moment angle manifolds and complexes}

We recall the notion of real moment angle complexes. Let $r$ be a positive integer and $K$ be a simplicial complex with vertex set $[r] = \lbrace 1,..., r \rbrace$. For each simplex $\sigma \in K$, we define 
$$ (D^1, S^0)^{\sigma}= \bigl \{ (x_1,..., x_r) \in (D^1)^r ~|~ x_i \in S^0 ~\mbox{when}~ i \notin \sigma \bigl \}. $$
Then the set 
$$ \mathbb{R}\mathcal{Z}_K:= \bigcup_{\sigma \in K}(D^1, S^0)^{\sigma} \subseteq (D^1)^r$$
is called the real moment angle complex of $K$. The space $\mathbb{R}\mathcal{Z}_K$ has a natural $\ZZ_2^r$-action induced from the $\ZZ_2^r$-action on $(D^1)^r$.

Let $P$ be a simple polytope with facets $\{ F_1,...,F_r\}$. Then the set 
$$ K_P:= \bigl \{ \sigma= \{ i_1,..., i_k\} ~|~ F_{i_1} \cap \cdots \cap F_{i_k} \neq \emptyset \bigl \}$$
is a simplicial complex on $\{ 1,..., r \}$, see \cite[Chapter 1]{BP02}. The set $K_P$ is called the dual of $P$, and $\mathbb{R}\mathcal{Z}_{K_P}$ has a manifold structure. The space $\mathbb{R}\mathcal{Z}_{K_P}$ is called the real moment angle manifold for $P$.
\begin{proposition} \label{small_cover_moment_manifold}
    Let $M^n(P, \lambda)$ be a small cover. Then there is a subgroup $\ZZ_{\lambda}$ of $\ZZ_2^r$ of rank $r-n$ such that $\ZZ_{\lambda}$ acts on $\mathbb{R}\mathcal{Z}_{K_P}$ freely and  $\mathbb{R}\mathcal{Z}_{K_P}/ \ZZ_{\lambda} \cong M^n(P, \lambda)$.
\end{proposition}
\begin{proof}
    This is similar to the proof of \cite[Proposition 6.5]{BP02}, and \cite[Proposition 2.4]{SZ22}.
\end{proof}
We note that $\mathbb{R}\mathcal{Z}_{K_{\Delta^n}}= \mathbb{S}^n$ and $\mathbb{R}\mathcal{Z}_{K_{P_1 \times P_2}}= \mathbb{R}\mathcal{Z}_{K_{P_1}} \times \mathbb{R}\mathcal{Z}_{K_{P_2}}$. Let $M^n(P, \lambda)$ be a small cover over an $n$-dimensional polytope $\prod_{j=1}^{m} \Delta^{n_j}$ for $j= 1,...,m$. Then, the number of facets of $\prod_{j=1}^{m} \Delta^{n_j}$ is $n+m$, and the real moment angle manifold $\mathbb{R}\mathcal{Z}_{K_P}$ is $\prod_{j=1}^{m} \mathbb{S}^{n_j}$. By Proposition \ref{small_cover_moment_manifold}, $M^n(P, \lambda)$ can be realized as the orbit space of the moment angle manifold $\prod_{j=1}^m \mathbb{S}^{n_j}$ by a free $\mathbb{Z}_2^m$-action. More precisely, the action of $\mathbb{Z}_2^m$ on $\prod_{j=1}^{m} \mathbb{S}^{n_j}$ is given by
\begin{align}
   & (g_1, g_2,...,g_m)((x_0^1,...,x_{n_1}^1),..., (x_0^m,..., x_{n_m}^m))  \\
    &= ((g_1x_0^1, (g_1^{a_{11}^1} \cdots g_m^{a_{m1}^1}) \cdot x_1^1,...,(g_1^{a_{1{n_1}}^1} \cdots g_m^{a_{m{n_1}}^1}) \cdot x_{n_1}^1),..., \nonumber \\
   & (g_m \cdot x_0^m, (g_1^{a_{11}^m} \cdots g_m^{a_{m1}^m}) \cdot x_1^m,...,(g_1^{a_{1n_m}^m} \cdots g_m^{a_{mn_m}^m}) \cdot x_{n_m}^m))\nonumber
\end{align}
\noindent where $(g_1, g_2,...,g_m) \in \mathbb{Z}_2^m$ and $(x_0^j,..., x_{n_j}^j) \in \mathbb{S}^{n_j}$ for $j=1,..., m$, see \cite[Remark 2.3]{DU19}. This $\mathbb{Z}_2^m$-action on $\prod_{j=1}^{m} \mathbb{S}^{n_j}$ is free and one has $\prod_{j=1}^{m} \mathbb{S}^{n_j}/\mathbb{Z}_2^m \cong M^n(P, \lambda)$.

\section{Equivariant LS-Category of small covers}\label{sec_eqlscat_small_cover}

In this section, we recall some basics of LS-category following \cite{CLOT03}. Then, we compute the LS-category and the equivariant LS-category of a small cover over a simple polytope. Next, we compute the LS-category of the real moment angle manifold for $r$-gon and the equivariant LS-category of the real moment angle complex.

 Let $G$ be a compact topological group acting continuously on a Hausdorff topological space $Y$. In this case, $Y$ is called a $G$-space. A subset $U$ of a $G$-space $Y$ is called $G$-invariant if $GU \subseteq U$. The homotopy $H \colon U \times I \rightarrow Y$ is called $G$-homotopy if for any $g \in G, y \in U$ and $t \in I$, we have $gH(y, t)=H(gy, t)$. A $G$-invariant open subset $U$ of $Y$ is called $G$-categorical if there exists an equivariant homotopy $H \colon U \times I \rightarrow Y$ such that $H_0$ is the inclusion, and $H_1 \colon U \rightarrow Y$ has the image in a single $G$-orbit. In particular, $U$ is called categorical if $G$ is trivial. Here we denote the orbit of an element $y \in Y$ by $\mathcal{O}(y)$.

\begin{definition} \label{def:ls_category}
The equivariant LS-category of a $G$-space $Y$, denoted by $\mbox{cat}_G(Y)$, is the least positive number of $G$-categorical invariant open sets required to cover $Y$. If no such covering exists, then $\mbox{cat}_G(Y)= \infty$.

In particular, if $G$ is trivial, then  $\mbox{cat}_G(Y)$ is called the LS-category of $Y$, denoted by  $\mbox{cat}(Y)$. 
\end{definition}


Let $Y$ be a space and $R$ be a commutative ring. The least integer $n$ such that all $(n + 1)$-fold cup products vanish in $H^*(Y; R)$ is called the cup-length of $Y$ with coefficients in $R$, denoted by $\mbox{cl}_R (Y)$. If no such $n$ exists, we write $\mbox{cl}_R (Y)= \infty $. The cup-length gives a lower bound for LS-category, as follows:

\begin{proposition} \label{cl_less_lscat}
 The cup-length of a topological space $Y$ is less than the LS-category of $Y$, i.e., $\mbox{cl}_R(Y)+ 1 \leq \mbox{cat}(Y)$, see \cite[Proposition 1.5]{CLOT03}.
\end{proposition}

\begin{proposition}  \label{lscat_less_dim}
 If $Y$ is a manifold, then $\mbox{cat}(Y) \leq \mbox{dim}(Y)+1$, see \cite[Theorem 1.7]{CLOT03}.
\end{proposition}



\begin{theorem} \label{small_cover_lscat}
 Let $M^n(P, \lambda)$ be an $n$-dimensional small cover. Then $\mbox{cat}(M^n(P, \lambda))= n+1$.
\end{theorem}

\begin{proof}
 
 Since $P$ is a simple polytope, at each vertex, exactly $n$ many facets intersect. Let $v$ be a vertex of $P$, and $v= F_{s_1} \cap \cdots \cap F_{s_n}$ where $ F_{s_1},..., F_{s_n}$ are unique $n$ facets of $P$. Let $m_v= \rho^{-1}(v)$ and $M_i= \rho^{-1}(F_{s_i})$ for $i=1,..., n$. Here the $\ZZ_2^n$-action on $M^n(P, \lambda)$ is locally standard. So, $m_v$ is a fixed point, and $M_1,..., M_n$ intersect to $m_v$ transversely. Therefore the Poincare dual of $M_i$ represents a non-zero cohomology class in $H^1(M^n(P, \lambda); \ZZ_2)$. So by definition of cup-length, $n \leq \mbox{cl}_{\ZZ_2}(M^n(P, \lambda))$. Therefore by Proposition \ref{cl_less_lscat}, $n+1 \leq \mbox{cat}(M^n(P, \lambda))$.
 Also by Proposition \ref{lscat_less_dim}, we have $\mbox{cat}(M^n(P, \lambda)) \leq \mbox{dim}(M^n(P, \lambda))+ 1= n+1$. Hence $\mbox{cat}(M^n(P, \lambda))= n+1$.
\end{proof}

  




We remark that the LS-category of small covers has been studied in \cite{ML21}. However, it is written in Chinese. So, we write a proof.

We recall a result from \cite{BS19}, which helps us to calculate the equivariant LS-category of a small cover over a simple polytope.

\begin{proposition} \cite[Theorem 3.3]{BS19} \label{orbit_class_lower_bound_equ_lscat}
Let $Y$ be a $G$ space and $\bigl \{[\mathcal{O}(y_i)] \bigl\}_{i \in \mathcal{A}}$ be the collection of all minimal orbit classes in $Y$. Let
$$ Y_i= \bigcup_{\mathcal{O}(y) \in [\mathcal{O}(y_i)]} \mathcal{O}(y).$$
Then 
$$ \# \mathcal{A} \leq \sum_{i \in \mathcal{A}} \mbox{cat}_G (Y_i) \leq \mbox{cat}_G (Y) $$
where $ \# \mathcal{A}$ is the cardinality of $\mathcal{A}$.
\end{proposition}

\begin{theorem}
 Let $M^n(P,\lambda)$ be an $n$-dimensional small cover over a simple polytope $P$ with $k$ vertices. Then $\mbox{cat}_{\mathbb{Z}_2^n}({M^n(P, \lambda)})= k$.
\end{theorem}
\begin{proof}
Let $M:= M^n(P,\lambda)$. We know that there is a bijection between the fixed point set $M^{\mathbb{Z}_2^n}$ and $V(P)$. Since the fixed points are isolated and minimal orbits, by Proposition \ref{orbit_class_lower_bound_equ_lscat}, we have $\mbox{cat}_{\mathbb{Z}_2^n}(M) \geq |V(P)|$. So, it is enough to show that for any $v \in M^{\mathbb{Z}_2^n}$, there is a $\mathbb{Z}_2^n$-categorical subset $X_v$ of $M$ such that $M= \bigcup_{v \in M^{\mathbb{Z}_2^n}} X_v$. Let $\rho \colon M \rightarrow P$ be the orbit map.
Now for $v \in M^{\mathbb{Z}_2^n}$, let 
$$C_v= \bigcup_{\rho(v) \notin F}F,~ U_v=P-C_v, ~\mbox{and}~ X_v=\rho^{-1}(U_v),$$
where $F$ is a face of $P$. Here $X_v$ is $\mathbb{Z}_2^n$-invariant subset of $M$. Since $U_v$ is a convex subset of $P$, it is contractible to $v$. So there exists a homotopy $h \colon U_v \times I \rightarrow P$ such that $h(x, 0)=x$ and $h(x, 1)=v$ for all $x \in U_v$ and preserves the face structure of $U_v \times I$. So, for any face $F$ of $U_v$, we have $h(x, t) \in F$ for $x \in F,~ t \in I$. Thus, by Proposition 1.8 of \cite{DJ91}, we can say $X_v \cong ( U_v \times \mathbb{Z}_2^n )/ \sim$. Therefore $h$ induces a homotopy 
$$ h \times Id \colon U_v \times I \times \mathbb{Z}_2^n \rightarrow P \times \mathbb{Z}_2^n  $$
defined by $(x, (r', t)) \mapsto (h(x, r'), t)$. Since for each face $F$ of $U_v$, we have
$$x \in F \Rightarrow h(x, r') \in F, ~\mbox{for all}~ r' \in I,$$
$h \times Id$ induces a homotopy
$H \colon X_v \times I \rightarrow M$
with 
$([x,t],r') \mapsto [h(x, r'), t]$. Since
$$gH([x,t],r')= g[h(x,r'), t]=[h(x, r'), gt]=H([x, gt],r')=H(g[x, t], r'),$$
the map $H$ is a $\mathbb{Z}_2^n$-homotopy. Also $H(x, 0)=x, H(x, 1)= \rho^{-1}(v)={v}, \text{ for all } x \in X_v$. Thus $X_v$ is $\mathbb{Z}_2^n$-categorical open invariant subset of $M$. Since $\{X_v ~|~ v \in V(P)\}$ covers $M$, therefore $\mbox{cat}_{\mathbb{Z}_2^n}(M)=|V(P)|=k$.
\end{proof}



\begin{proposition}
 Let $P$ be an $r$-gon and $\mathbb{R}\mathcal{Z}_{K_P}$ be a moment angle manifold. Then ${\mbox{cat}(\mathbb{R}\mathcal{Z}_{K_P})= 3}$.
\end{proposition}
\begin{proof}
We know the cohomology ring $H^*(\mathbb{R}\mathcal{Z}_{K_P}; \ZZ_2)$ is generated by elements of degree only $0, 1$ and $2$. We can get two elements of degree $1$ such that their cup product is non-zero in $H^*(\mathbb{R}\mathcal{Z}_{K_P}; \ZZ_2)$, see \cite[Section 3]{Cai17}. Therefore $2 \leq \mbox{cl}_{\ZZ_2}(\mathbb{R}\mathcal{Z}_{K_P})$. Then, we have $3 \leq \mbox{cat}(\mathbb{R}\mathcal{Z}_{K_P})$. Also $\mbox{dim}(\mathbb{R}\mathcal{Z}_{K_P})= \mbox{dim}(P)=2$. Therefore, $\mbox{cat}(\mathbb{R}\mathcal{Z}_{K_P}) \leq 3$. Hence $\mbox{cat}(\mathbb{R}\mathcal{Z}_{K_P})= 3$.
\end{proof}
\begin{remark}
Let $K$ be a triangulated $d$-sphere for $d \leq 2$ or a connected sum of joins of such spheres. If $K$ is $k$-Golod over $\ZZ_2$ (i.e., length $k+1$ cup products of positive degree elements in $H^*(\mathbb{R}\mathcal{Z}_K; \ZZ_2)$ vanish), then $k \leq \mbox{cl}_{\ZZ_2}(\mathbb{R}\mathcal{Z}_K)$. Thus $k+1 \leq \mbox{cat}(\mathbb{R}\mathcal{Z}_K)$, see \cite[Theorem 4.2]{BG21}.
\end{remark}

\begin{theorem}
Let $S$ be the set of all maximal simplices of a simplicial complex $K$ on $[r]$. Then 
$$ { \mbox{cat}_{\mathbb{Z}_2^r}(\mathbb{R}\mathcal{Z}_K) = |S|}.$$
\end{theorem}

\begin{proof}
Note that if $\tau$ is a face of $\sigma$ in $K$, then $(D^1, S^0)^{\tau} \subseteq (D^1, S^0)^{\sigma}$. So we have
$$ \mathbb{R}\mathcal{Z}_K= \bigcup_{\sigma \in S}(D^1, S^0)^{\sigma} \subseteq (D^1)^r.$$
The topology on $\mathbb{R}\mathcal{Z}_K$ is the subspace topology of $(D^1)^r$. Also, any simplex of $K$ is a
face of a maximal simplex. So the set
	
$$ \lbrace (D^1, S^0)^{\sigma} ~|~ \sigma \in S \rbrace $$
is an open covering for $\mathbb{R}\mathcal{Z}_K$. Moreover, $(D^1, S^0)^{\sigma}$ is a $\mathbb{Z}_2^r$-invariant subset which is equivariantly contractible to the orbit $(S^0)^{\sigma}$ in $\mathbb{R}\mathcal{Z}_K$ where
$$ (S^0)^{\sigma}= \lbrace (x_1,..., x_r) \in \mathbb{R}\mathcal{Z}_K ~|~ x_i= 0~ \mbox{if}~ i \in \sigma ~\mbox{and}~ |x_i|= 1 ~ \mbox{if}~ i \notin \sigma \rbrace. $$
So we obtain that
$$ \mbox{cat}_{\mathbb{Z}_2^r}(\mathbb{R}\mathcal{Z}_K) \leq |S|.$$
Note that the set $ \lbrace (S^0)^{\sigma} ~|~ \sigma \in S \rbrace$ is the set of all minimal orbits of $\mathbb{R}\mathcal{Z}_K$ with respect to $\mathbb{Z}_2^r$-action. So, by Proposition \ref{orbit_class_lower_bound_equ_lscat}, we have
$$\mbox{cat}_{\mathbb{Z}_2^r}(\mathbb{R}\mathcal{Z}_K) \geq |S|.$$
\end{proof}

\section{Topological complexity of small covers} \label{sec_higher_tc_small_cover}

 In this section, we recall the definition of topological complexity and zero-divisors-cup-length. Next, we try to give bounds for the topological complexity of a class of small covers over a product of simplices and real Bott manifolds. 
 
 \begin{definition} 
 \label{TC_definition}
Let $Y$ be a path-connected space. The topological complexity of the motion planning in $Y$ is the least integer $k$ such that $Y \times Y$ can be covered by $k$ open subsets $U_1,..., U_k$ on each of which there exists a section $s_i \colon U_i \rightarrow PY$ such that $\pi \circ s_i$ is homotopic to the inclusion $id_{U_i}$. If no such integer exists, then we set $\mathbf{TC}(Y)= \infty$.
\end{definition}

We note that in the above definition, we consider non-normalized topological complexity.
The cup product map
\begin{equation}\label{tensor_map}
\cup \colon H^*(Y; R) \otimes H^*(Y; R) \rightarrow H^*(Y; R)
\end{equation}
is an algebra homomorphism whose kernel is called the ideal of zero-divisors of $H^*(Y; R)$. The multiplicative structure on the left in (\ref{tensor_map}) is given by the formula ${(\alpha \otimes \beta) \cdot (\gamma \otimes \delta)= (-1)^{|\beta|.|\gamma|} \alpha \gamma \otimes \beta \delta}$. Here $|\beta|$ and $|\gamma|$ denote the degrees of the cohomology classes $\beta$ and $\gamma$ respectively.

\begin{definition}
The zero-divisors-cup-length of $H^*(Y; R)$, denoted by $\mbox{zcl}_R(Y)$, is the length of the longest nontrivial product in the ideal of the zero-divisors of $H^*(Y; R)$.
\end{definition}
The following proposition gives a lower bound and an upper bound for $\mathbf{TC}(Y)$.

\begin{proposition} \cite[Theorem 4,5,7]{Far03} \label{zcl_less_tc_less_dim}
If $Y$ is a manifold, then we have
$$ \mbox{max} \lbrace \mbox{cat}(Y), \mbox{zcl}_R(Y)+1 \rbrace \leq \mathbf{TC}(Y) \leq 2 \mbox{dim}(Y) +1.$$
\end{proposition}

\begin{proposition} \cite[Theorem 11]{Far03} \label{tc_product_inequality}
For any path-connected metric spaces  $Y_1, ..., Y_m$, we have
$$\mathbf{TC}(Y_1 \times \cdots \times Y_m) \leq \mathbf{TC}(Y_1) + \cdots + \mathbf{TC}(Y_m)- (m-1).$$
\end{proposition} 

\begin{proposition} \label{tc_small_cover_upper_lower}
 Let $M^n(P, \lambda)$ be an $n$-dimensional small cover. Then $$ n+1 \leq \mathbf{TC}(M^n(P, \lambda)) \leq 2n+1.$$
\end{proposition}

\begin{proof}
The proof follows from Theorem \ref{small_cover_lscat} and Proposition \ref{zcl_less_tc_less_dim}.
\end{proof}

Now we calculate the topological complexity of several small covers over a product of two simplices. Consider the set $\mathcal{S}:= \{ n \in \mathbb{N} ~|~ \binom{n}{i} ~\mbox{is even for}~ 0<i<n \}$. Let $n \leq 2^r-1<2n$, ${n_j \leq 2^{r_j}-1<2n_j}$ for $j=1,2$, and $n=n_1+n_2$. Note that $\binom{2^s}{i}$ is even for $0 < i < 2^s$.
\begin{theorem}\label{tc_small_cover_special}
Let $M^n(P, \lambda)$ be a small cover other than $\mathbb{RP}^{n_1} \times \mathbb{RP}^{n_2}$ over $P= \Delta^{n_1} \times \Delta^{n_2}$.
\begin{enumerate}
 \item \label{tc_small_cover_special_case1} Let $n_2 \in \mathcal{S}$ with $n_2 > n_1$. Then $2^{r_1}+2^{r_2}-1 \leq \mathbf{TC}(M^n(P, \lambda))$.
 \item \label{tc_small_cover_special_case2} Let $n_2 \in \mathcal{S}$ with $n_2$ divides $n_1$. Then $2^r \leq \mathbf{TC}(M^n(P, \lambda))$.
\item \label{tc_small_cover_special_case3} Let $n_2 \in \mathcal{S}+1$ with $n_2 > n_1+1$. Then, $ 2^r \leq \mathbf{TC}(M^n(P, \lambda))$.

\item \label{tc_small_cover_special_case4} Let $n_2 \in \mathcal{S}+2$ with $n_2 > n_1+2$. Then, $2^r \leq \mathbf{TC}(M^n(P, \lambda))$.
\end{enumerate}
In particular, if $n=2^{s-1}$, then for the cases \eqref{tc_small_cover_special_case2}, \eqref{tc_small_cover_special_case3}, and \eqref{tc_small_cover_special_case4}, we have $$2n \leq \mathbf{TC}(M^n(P, \lambda)) \leq 2n+1.$$
\end{theorem}
\begin{proof}
 In the cohomology ring $H^*(M^n(P, \lambda); \ZZ_2)$ described in Proposition \ref{cohomoring_smallcover_polyt}, the ideal $\Tilde{I}$ is generated by $\alpha_1=x_1x_2 \cdots x_{n_1}y_1$ and $\alpha_2=x_{n_1+1}x_{n_1+2} \cdots x_{n_1+n_2}y_2$. The ideal $\Tilde{J}$ is generated by 
$$x_1= x_2= \cdots =x_{n_1}=y_1,
~\mbox{and}~
x_{n_1+1}=x_{n_1+2}= \cdots= x_{n_1+n_2}= y_1+ y_2.$$

Therefore $\alpha_1= y_1^{n_1+1}= 0$ and $y_1^{n_1} \neq 0$ in $H^*(M^n(P, \lambda); \ZZ_2)$. 
From \eqref{at_v00000}, and the Poincare duality, we have ${x_1  \cdots x_{n_1}x_{n_1+1}  \cdots x_{n_1+n_2} \neq 0}$. So, $y_1^{n_1}(y_1+y_2)^{n_2} \neq 0$. Now,
 \begin{align*}
y_1^{n_1}(y_1+y_2)^{n_2} &=  y_1^{n_1} \{ y_2^{n_2}+y_1 \cdot f(y_1, y_2)\} ~(\mbox{where}~ f(y_1, y_2) ~\mbox{is a function of}~ y_1~ \mbox{and}~ y_2) \\
&= y_1^{n_1} y_2^{n_2}+y_1^{n_1+1} \cdot f(y_1, y_2) \\
&=  y_1^{n_1} y_2^{n_2} ~( \mbox{as}~ y_1^{n_1+1}= 0).
\end{align*}
Therefore we get the following:
\begin{equation} \label{y1n1y2n2_0}
 y_1^{n_1} y_2^{n_2} \neq 0.  
\end{equation}

 Let $\mathfrak{a}_j: =1 \otimes y_j- y_j \otimes 1$ for $j=1, 2$. Then $\mathfrak{a}_j$ is in the ideal of the zero-divisors of $H^*(M^n(P, \lambda); \ZZ_2)$. Let $c=2^r-1$, and $c_j= 2^{r_j}-1$ for $j=1, 2$.
\begin{enumerate}
\item \label{TC_main_proof_case1}
Let $n_2 \in \mathcal{S}$ with $n_2 > n_1$. Here,
\begin{align*}
\alpha_2 &=  (y_1+y_2)^{n_2}y_2\\
&= (y_1^{n_2}+y_2^{n_2})y_2 ~(\mbox{as}~ \binom{n_2}{i} ~\mbox{is even for}~ 0<i<n) \nonumber\\
&= y_2^{n_2+1} ~( \text{ as } y_1^{n_2}=0, ~ \mbox{since}~ n_2 \geq n_1+1 ~\mbox{and}~ y_1^{n_1+1}= 0).\nonumber
\end{align*}
Therefore $y_2^{n_2+1}=0$. Now, for $j=1, 2$,
$$\mathfrak{a}_j^{c_j}=  (1 \otimes y_j- y_j \otimes 1) ^{c_j}= \sum_{k_j=0}^{c_j} (-1)^{c_j-k_j} \binom{c_j}{k_j} (y_j^{c_j-k_j} \otimes y_j^{k_j}).$$

Now, the binomial coefficient $\binom {2^{r_j}-1} {i_j}$ is odd for all $0 \leq i_j \leq c_j$ for $j=1, 2$. The binomial expansion of  $\mathfrak{a}_j^{c_j}$ contains the term $(y_j^{c_j-{n_j}} \otimes y_j^{n_j})$ which is non-zero. Now, by \eqref{y1n1y2n2_0}, $y_1^{n_1}y_2^{n_2} \neq 0$. So, $\mathfrak{a}_1^{c_1} \mathfrak{a}_2^{c_2}$ contains the term $y_1^{n_1}y_2^{n_2} \otimes (y_1^{c_1-{n_1}} y_2^{c_2-{n_2}})$ which is non-zero and there is no other term of this form in the expression of $\mathfrak{a}_1^{c_1} \mathfrak{a}_2^{c_2}$.

Hence zero-divisors-cup-length of $H^*(M^n(P, \lambda); \ZZ_2)$ is greater than or equal to $c_1+c_2$. Therefore, by Proposition \ref{zcl_less_tc_less_dim}, we have, $2^{r_1}+2^{r_2}-1 \leq \mathbf{TC}(M^n(P, \lambda))$.

\item \label{TC_main_proof_case2} 
Now, consider the case when $n_2 \in \mathcal{S}$ with $n_2$ divides $n_1$. Let $n_1= \Bar{n}n_2$ for some $\Bar{n} \in \ZZ$. Since $\binom{n_2}{i}$ is even for $0<i<n_2$, so,
$$ \alpha_2 = (y_1+y_2)^{n_2}y_2= (y_1^{n_2}+y_2^{n_2})y_2= y_1^{n_2}y_2+ y_2^{n_2+1}.$$
So, $ y_2^{n_2+1}= y_1^{n_2}y_2$ in $H^*(M^n(P, \lambda); \ZZ_2)$. Thus,
\begin{align*}
y_1^{n_1}y_2^{n_2}= y_1^{\Bar{n}n_2}y_2^{n_2}= y_1^{\Bar{n}n_2}y_2^{\Bar{n}}y_2^{n_2-\Bar{n}}= y_2^{\Bar{n}n_2+\Bar{n}}y_2^{n_2-\Bar{n}}= y_2^{\Bar{n}n_2+n_2}=y_2^{n_1+n_2}=y_2^n. \nonumber
\end{align*}
Now, by \eqref{y1n1y2n2_0}, $y_1^{n_1}y_2^{n_2} \neq 0$. Thus, $y_2^n \neq 0$.
Now,
$$\mathfrak{a}_2^c= (1 \otimes y_2- y_2 \otimes 1) ^c= \sum_{k=0}^{c} (-1)^{c-k} \binom{c}{k} (y_2^{c-k} \otimes y_2^k).$$
Therefore, by similar arguments as in \eqref{TC_main_proof_case1}, we get $2^r \leq \mathbf{TC}(M^n(P, \lambda))$.


\item  \label{TC_main_proof_case3}
Let $n_2 \in \mathcal{S}+1$ with $n_2 > n_1+1$. Now,
\begin{align*}
\alpha_2 &=  (y_1+y_2)^{n_2}y_2\\
&=  (y_1+y_2)^{n_2-1}(y_1+y_2)y_2 \\
&= (y_1^{n_2-1}+y_2^{n_2-1})(y_1+y_2)y_2 ~(\mbox{as}~ \binom{n_2-1}{i} ~\mbox{is even for}~ 0<i<n_2-1) \nonumber\\
&= (y_1+y_2)y_2^{n_2} ~( \text{ as } y_1^{n_2-1}=0, ~ \mbox{since}~ n_2-1 \geq n_1+1 ~\mbox{and}~ y_1^{n_1+1}= 0)\nonumber \\
&= y_1y_2^{n_2}+ y_2^{n_2+1}.  \nonumber
\end{align*}
So, $y_2^{n_2+1}= y_1y_2^{n_2}$ in $H^*(M^n(P, \lambda); \ZZ_2)$. Thus,
\begin{align*}
y_1^{n_1}y_2^{n_2}= y_1^{n_1-1}(y_1y_2^{n_2})= y_1^{n_1-1}y_2^{n_2+1}= \cdots= y_2^{n_1+n_2}= y_2^n. \nonumber
\end{align*}
Therefore, by similar arguments as in \eqref{TC_main_proof_case2}, we get $2^r \leq \mathbf{TC}(M^n(P, \lambda))$.


\item  \label{TC_main_proof_case4}
Let $n_2 \in \mathcal{S}+2$ with $n_2 > n_1+2$. Now,
\begin{align*}
\alpha_2 &=  (y_1+y_2)^{n_2}y_2\\
&=  (y_1+y_2)^{n_2-2}(y_1+y_2)^2y_2 \\
&= (y_1^{n_2-2}+y_2^{n_2-2})(y_1^2+y_2^2)y_2 ~(\mbox{as}~ \binom{n_2-2}{i} ~\mbox{is even for}~ 0<i<n_2-2) \nonumber\\
&= (y_1^2+y_2^2)y_2^{n_2-1} ~( \text{ as } y_1^{n_2-2}=0, ~ \mbox{since}~ n_2-2 \geq n_1+1 ~\mbox{and}~ y_1^{n_1+1}= 0)\nonumber \\
&= y_1^2y_2^{n_2-1}+ y_2^{n_2+1}.  \nonumber
\end{align*}
So, $y_2^{n_2+1}= y_1^2y_2^{n_2-1}$ in $H^*(M^n(P, \lambda); \ZZ_2)$. Let $n_1$ be even. Then,
\begin{align*}
y_1^{n_1}y_2^{n_2}= y_1^{n_1-2}(y_1^2y_2^{n_2-1})y_2= y_1^{n_1-2}y_2^{n_2+1}y_2= \cdots= y_2^{n_1+n_2}= y_2^n. \nonumber
\end{align*}
Therefore, by similar arguments as in \eqref{TC_main_proof_case2}, we get $2^r \leq \mathbf{TC}(M^n(P, \lambda))$.

Now, let $n_1$ be odd. Then,
\begin{align*}
y_1^{n_1}y_2^{n_2}= y_1^{n_1-2}(y_1^2y_2^{n_2-1})y_2= y_1^{n_1-2}y_2^{n_2+1}y_2= \cdots= y_1y_2^{n_1+n_2-1}= y_1y_2^{n-1}. \nonumber
\end{align*}
Now, by \eqref{y1n1y2n2_0}, $y_1^{n_1}y_2^{n_2} \neq 0$. Therefore $y_1y_2^{n-1} \neq 0$. So, $y_2^{n-1} \neq 0$. We know $n \leq 2^r-1 < 2n$, i.e., $2^r-1 \leq 2n-1$. If $2^r-1= 2n-1$, then $n$ is even. Since $n_1$ is odd, so $n_2$ is odd, which is not true. Therefore,  $2^r-1 \leq 2n-2$. 
Now,
$$\mathfrak{a}_2^c= (1 \otimes y_2- y_2 \otimes 1) ^c= \sum_{k=0}^{c} (-1)^{c-k} \binom{c}{k} (y_2^{c-k} \otimes y_2^k).$$

Here, the binomial expansion of  $\mathfrak{a}_2^c$ contains the term $(y_2^{c-n+1} \otimes y_2^{n-1})$ which is non-zero, and there is no other same term in the expression of $\mathfrak{a}_2^c$. So $\mathfrak{a}_2^c$ is non-zero. Therefore, by similar arguments as in \eqref{TC_main_proof_case2}, we get $2^r \leq \mathbf{TC}(M^n(P, \lambda))$.

In particular, if $n=2^{s-1}$ then $n \leq 2^{s}-1<2n$. Therefore, by \eqref{TC_main_proof_case2}, \eqref{TC_main_proof_case3}, and \eqref{TC_main_proof_case4} we have ${2^s \leq \mathbf{TC}(M^n(P, \lambda))}$. Thus $ 2n \leq \mathbf{TC}(M^n(P, \lambda)) \leq 2n+1$ for any $s \geq 1$.

\end{enumerate}
\end{proof}

\begin{proposition} \cite[Corollary 8.1, 8.2]{FTY03} \label{tc_137_powerof2}
 If $n$ equals $1,3$ or $7$, then $\mathbf{TC}(\mathbb{RP}^n)= n+1$, and if $n$ is a power of $2$, then $\mathbf{TC}(\mathbb{RP}^n)= 2n$. 
\end{proposition}
\begin{proposition} \label{tc_small_cover_special_case}
Let the small cover $M^n(P, \lambda)$ over the polytope $P= \prod_{j=1}^m \Delta^{n_j}$ is of the form $\mathbb{RP}^{n_1} \times \cdots \times \mathbb{RP}^{n_m}$.
\begin{enumerate}
 \item If $n_j= 2^{s_j-1}$ for all $j \in \{1,2,\dots, m\}$, then $\mathbf{TC}(M^n(P, \lambda))= 2^{s_1}+ \cdots + 2^{s_m}- (m-1)$.

\item If $n_j= 1,3$ or $7$ for all $j \in \{1,2, \dots, m\}$, then $\mathbf{TC}(M^n(P, \lambda))= n+1$.
\end{enumerate}
\end{proposition}
\begin{proof}

Here $M^n(P, \lambda)= \mathbb{RP}^{n_1} \times \cdots \times \mathbb{RP}^{n_m}$. In the cohomology ring $H^*(M^n(P, \lambda); \ZZ_2)$, the ideal $\Tilde{I}$ is generated by $\alpha_j= x_{\mathcal{N}_{j-1}+1} x_{\mathcal{N}_{j-1}+2} \cdots x_{\mathcal{N}_j} y_j $ for $j=1,...,m$, and the ideal $\Tilde{J}$ is generated by 
\begin{equation} \label{y1n1ymnm_neq_0}
x_{\mathcal{N}_{j-1}+1}= x_{\mathcal{N}_{j-1}+2}= \cdots= x_{\mathcal{N}_j}= y_j.
\end{equation}

Therefore, $\alpha_j= y_j^{n_j+1}= 0$ in $H^*(M^n(P, \lambda); \ZZ_2)$ for $j=1,...,m$. Let $$\mathfrak{a}_j: =1 \otimes y_j- y_j \otimes 1$$ for $j=1,..., m$. Then $\mathfrak{a}_j$ belongs to the ideal of the zero-divisors of $H^*(M^n(P, \lambda); \ZZ_2)$. Let $c_j=2^{s_j}-1$ for $j=1,...,m$. 
Now, by Lemma \ref{yj_nj_neq_0}, $y_j^{n_j} \neq 0$. From \eqref{at_v00000}, and the Poincare duality, we have $${x_1  \cdots x_{\mathcal{N}_1}x_{\mathcal{N}_1+1}  \cdots x_{\mathcal{N}_2}x_{\mathcal{N}_2+1} \cdots x_{\mathcal{N}_{m-1}+1} \cdots x_{\mathcal{N}_m} \neq 0}.$$ Using \eqref{y1n1ymnm_neq_0}, we have, $y_1^{n_1} \cdots y_m^{n_m} \neq 0$. 
Therefore, by similar arguments as in Theorem \ref{tc_small_cover_special}\eqref{TC_main_proof_case1}, we have $\mathfrak{a}_1^{c_1} \mathfrak{a}_2^{c_2} \cdots \mathfrak{a}_m^{c_m} \neq 0$. Thus, $ c_1+ \cdots + c_m+ 1 \leq \mathbf{TC}(M^n(P, \lambda))$. That is $$ 2^{s_1}+ \cdots + 2^{s_m}- (m-1) \leq \mathbf{TC}(M^n(P, \lambda)).$$

Also, from Proposition \ref{tc_product_inequality}, we have
 \begin{equation} \label{tc_rpn_inequality}
 \mathbf{TC}(M^n(P, \lambda)) \leq \mathbf{TC}(\mathbb{RP}^{n_1})+ \cdots + \mathbf{TC}(\mathbb{RP}^{n_m})- (m-1).
 \end{equation}
\begin{enumerate}
\item If $n_j=2^{s_j-1}$ then $n_j \leq 2^{s_j}-1<2n_j$. Now by Proposition \ref{tc_137_powerof2}, we have $\mathbf{TC}(\mathbb{RP}^{n_j})= 2^{s_j}$. So, the right inequality can be obtained using \eqref{tc_rpn_inequality}. Hence $$\mathbf{TC}(M^n(P, \lambda))= 2^{s_1}+ \cdots + 2^{s_m}- (m-1).$$

\item If $n_j= 1,3$ or $7$, then there exists some $s_j$ which satisfies $n_j \leq 2^{s_j}-1 < 2n_j$ and $n_j+1= 2^{s_j}$. So, $2^{s_1}+ \cdots + 2^{s_m}= n_1+ \cdots+ n_m+m= n+m$. 
Thus, we have $n+m-(m-1)=n+1 \leq \mathbf{TC}(M^n(P, \lambda))$.  By Proposition \ref{tc_137_powerof2}, we have $\mathbf{TC}(\mathbb{RP}^{n_j})= n_j+1$. So, the right inequality can be obtained using \eqref{tc_rpn_inequality}. Hence $$\mathbf{TC}(M^n(P, \lambda))= n+1.$$
\end{enumerate}
\end{proof}
We remark that if $M^n(P, \lambda)$ is not $\mathbb{RP}^{n_1} \times \cdots \times \mathbb{RP}^{n_m}$ then computation of $\mathbf{TC}(M^n(P, \lambda))$ is a challenging problem.


    
We recall that the $n$-th stage real Bott manifold is a small cover $M^n(P, \lambda)$ over the polytope $P= \prod_{j=1}^n \Delta^{n_j}$ where $n_j=1$ for $j=1,...,n$ and $\lambda$ be as in \eqref{Eq_define_lambda_on_P}. In this case, the elements of the $(n \times n)$ matrix coming from \eqref{lower_triangular_matrix} are scalars. Note that the diagonal elements of this matrix are 1 follows from the definition of $\lambda$. Since the Bott matrix is unique up to conjugation, different ${\beta_l^m}'s$ give different real Bott manifolds up to equivariant diffeomorphism. Now we calculate some lower bounds (possibly tight) of the topological complexity of the real Bott manifolds.

\begin{theorem}\label{real_bott_tc_matrix_change}
For $n \geq 3$, let the elements $\beta_k^{k+1}$ in the Bott matrix \eqref{lower_triangular_matrix} be $1$ for $k=1,..., n-1$, and the remaining elements $\beta_l^m$ be zero for $l= 1,...,n-2$, and $m= 3,...,n$. If $n \leq 2^r-1 < 2n$, then the topological complexity of the real Bott manifold $M^n(P, \lambda)$ is greater than or equal to $2^r$. In particular, if $n= 2^{s-1}$, then
$2n \leq \mathbf{TC}(M^n(P, \lambda)) \leq 2n+1$.
\end{theorem}
\begin{proof}
   In the cohomology ring $H^*(M^n(P, \lambda); \ZZ_2)$ as in Proposition \ref{cohomoring_smallcover_polyt}, the generators of the ideal $\Tilde{I}$ are $\alpha_j= x_jy_j$ for $j=1,...,n$, and the ideal $\Tilde{J}$ is generated by the elements $x_1+ y_1$ and $x_j+ y_j+ y_{j-1}$ for $j= 2,...,n$. Now for  $j \in \{ 2,...,n \}$,
$$\alpha_j=  x_jy_j= (y_{j-1}+ y_j)y_j= y_{j-1}y_j+ y_j^2. $$

Our claim is that $y_j^2 \neq 0$ for $j= 2,...,n$. For this, it is enough to show that $y_{j-1}y_j \neq 0$ for $j= 2,...,n$ as $\alpha_j=0$ in $H^*(M^n(P, \lambda); \ZZ_2)$. Note that in this case, $P= \prod_1^n \Delta^1$, an $n$-cube. So, the facets corresponding to the indeterminates $x_j$ and $y_j$ don't intersect. But the facets corresponding to the indeterminates $y_{j-1}$ and $y_j$ intersect to an $(n-2)$-dimensional face. So $y_{j-1}y_j$ is non-zero in $H^*(M^n(P, \lambda); \ZZ_2)$. Therefore $y_j^2 \neq 0$ for $j=2,...,n$.

Since $P$ is an $n$-dimensional simple polytope, there is a vertex where the facets corresponding to the indeterminates $y_1, y_2,..., y_n$ intersect. In other words, $y_1 y_2 \cdots y_n \neq 0$, by Poincare duality. From the relation $y_{j-1}y_j= y_j^2$ for $j= 2,...,n$, we have $y_1 y_2 \cdots y_n= y_n^n$. Therefore, $y_n^n \neq 0$.

Let $\mathfrak{a}_n:= 1 \otimes y_n- y_n \otimes 1$. Then $\mathfrak{a}_n$ is in the ideal of the zero-divisors of $H^*(M^n(P, \lambda); \ZZ_2)$. Then,
\begin{align*}
\mathfrak{a}_n^{2^r-1} &=  (1 \otimes y_n- y_n \otimes 1) ^{2^r-1}\nonumber \\
&= \sum_{k=0}^{2^r-1} (-1)^{2^r-1-k} \binom{2^r-1}{k} (1 \otimes y_n)^k (y_n \otimes 1)^{2^r-1-k} \nonumber\\
&=  \sum_{k=0}^{2^r-1} (-1)^{2^r-1-k} \binom{2^r-1}{k} (y_n^{2^r-1-k} \otimes y_n^k).\nonumber
\end{align*}

The binomial coefficients $\binom {2^r-1} {i}$ are odd for all $0 \leq i \leq 2^r-1$. The binomial expansion of  $\mathfrak{a}_n^{2^r-1}$ contains the term $(y_n^{2^r-1-n} \otimes y_n^{n})$ which is non-zero and there is no other term of this form in the expression of $\mathfrak{a}_n^{2^r-1}$. So $\mathfrak{a}_n^{2^r-1}$ is nonzero. Therefore zero-divisors-cup-length of $H^*(M^n(P, \lambda); \ZZ_2)$ is greater than or equal to $2^r-1$. Hence, by Proposition \ref{zcl_less_tc_less_dim}, we have $2^r \leq \mathbf{TC}(M^n(P, \lambda)) $. 

If $n= 2^{s-1}$ then $r$ satisfies $n \leq 2^s-1 < 2n$. Thus, $2^s=2n \leq \mathbf{TC}(M^n(P, \lambda)) \leq 2n+1$.
\end{proof}

We recall that for $n=3$, the Bott matrix is given by
$\begin{pmatrix}
1 & 0 & 0 \\
\beta_1^2 & 1 & 0 \\
\beta_1^3 & \beta_2^3 & 1
\end{pmatrix}$. We denote the corresponding real Bott manifold $M^3(P, \lambda)$ by $M^3(\beta_1^2, \beta_1^3, \beta_2^3)$.

\begin{theorem}
${5 \leq \mathbf{TC}(M^3(1, 0, 0)), \mathbf{TC}(M^3(0, 1, 0)), \mathbf{TC}(M^3(0, 0, 1)), \mathbf{TC}(M^3(0, 1, 1)) \leq 7}$.
\end{theorem}
\begin{proof}
The generators of the ideal $\Tilde{I}$ in Proposition \ref{cohomoring_smallcover_polyt} are $\alpha_j= x_jy_j$ where $j=1, 2, 3$. Let $$\mathfrak{a}_j:= 1 \otimes y_j- y_j \otimes 1$$ for $j= 1, 2, 3$. Then $\mathfrak{a_j}$ is in the ideal of the zero-divisors of $H^*(M^3(\beta_1^2, \beta_1^3, \beta_2^3); \ZZ_2)$. Note that $\mathbf{TC}(M^3(\beta_1^2, \beta_1^3, \beta_2^3)) \leq 7$ by Proposition \ref{tc_small_cover_upper_lower}. The manifolds $M^3(1, 0, 0), M^3(0, 1, 0), M^3(0, 0, 1)$, and $ M^3(0, 1, 1)$ are diffeomorphic to each other by \cite[Theorem 4]{Naz11}.

Consider the real Bott manifold $M^3(1, 0, 0)$. Then from Proposition \ref{cohomoring_smallcover_polyt}, the ideal $\Tilde{J}$ is generated by the elements $x_1+ y_1, x_2+ y_1+y_2$, and $x_3+ y_3$. So, $x_1= y_1$, $x_2= y_1+y_2$, and $x_3= y_3$ in $H^*(M^3(1, 0, 0); \ZZ_2)$. Therefore, we have $y_1^2= y_3^2= 0$, and $y_2^2= y_1y_2$. 
Now,
\begin{align*}
\mathfrak{a}_2^3 \mathfrak{a}_3 &=  (1 \otimes y_2^3- y_2 \otimes y_2^2+ y_2^2 \otimes y_2- y_2^3 \otimes 1) (1 \otimes y_3- y_3 \otimes 1)  \nonumber \\
&=  y_1y_2 \otimes y_2y_3+ y_2y_3 \otimes y_1y_2- y_1y_2y_3 \otimes y_2-  y_2 \otimes y_1y_2y_3.\nonumber
\end{align*}
So, the product $\mathfrak{a}_2^3 \mathfrak{a}_3$ contains an element $y_1y_2 \otimes y_2y_3$ which is non-zero. Therefore, the zero-divisors-cup-length of $H^*(\mathbf{TC}(M^3(1, 0, 0)); \ZZ_2)$ is greater than or equal to $4$. Hence by Proposition \ref{zcl_less_tc_less_dim}, we have $5 \leq \mathbf{TC}(M^3(1, 0, 0))$.

\end{proof}

We remark that $M^3(1, 1, 0)$ is the $3$-dimensional Klein Bottle, and \cite[Theorem 3.1]{DS23} gives $\mathbf{TC}(M^3(1, 1, 0))= 6$.

\begin{theorem} \label{bott_mfd_36}
$6 \leq \mathbf{TC}(M^3(1, 0, 1)), \mathbf{TC}(M^3(1, 1, 1)) \leq 7$.
\end{theorem}

\begin{proof}
The generators of the ideal $\Tilde{I}$ in Proposition \ref{cohomoring_smallcover_polyt} are $\alpha_j= x_jy_j$ where $j=1, 2, 3$. Let $\mathfrak{a}_j:= 1 \otimes y_j- y_j \otimes 1$ for $j= 1, 2, 3$. Then $\mathfrak{a_j}$ is in the ideal of the zero-divisors of $H^*(M^3(P, \lambda); \ZZ_2)$. Note that $\mathbf{TC}(M^3(\beta_1^2, \beta_1^3, \beta_2^3)) \leq 7$ by Proposition \ref{tc_small_cover_upper_lower}. The manifolds $M^3(1, 0, 1)$ and $M^3(1, 1, 1)$ are diffeomorphic by \cite[Theorem 4]{Naz11}.


Consider the real Bott manifold $M^3(1, 0, 1)$. Then from Proposition \ref{cohomoring_smallcover_polyt}, the ideal $\Tilde{J}$ is generated by the elements $x_1+ y_1, x_2+ y_1+ y_2$, and $x_3+ y_2+ y_3$. So, $x_1= y_1, x_2= y_1+ y_2$, and $x_3= y_2+ y_3$ in $H^*(M^3(1, 0, 1); \ZZ_2)$.
Now,
\begin{align*}
\mathfrak{a}_2^2 \mathfrak{a}_3^3 &=  (1 \otimes y_2^2+ y_2 \otimes 1) (1 \otimes y_3^3- y_3 \otimes y_3^2+ y_3^2 \otimes y_3- y_3^3 \otimes 1)  \nonumber \\
&= (y_1y_2+ y_2y_3) \otimes y_1y_2y_3- y_1y_2y_3 \otimes (y_1y_2+ y_2y_3).\nonumber
\end{align*}
So, the product $\mathfrak{a}_2^2 \mathfrak{a}_3^3$ contains an element $y_1y_2 \otimes y_1y_2y_3$ which is non-zero. Therefore the zero-divisors-cup-length of $H^*(M^3(1, 0, 1); \ZZ_2)$ is greater than or equal to $5$. Hence by Proposition \ref{zcl_less_tc_less_dim}, we have $6 \leq \mathbf{TC}(M^3(1, 0, 1))$. 



\end{proof}

Now, for $n=4$, the Bott matrix is given by
$\begin{pmatrix}
1 & 0 & 0 & 0 \\
\beta_1^2 & 1 & 0 & 0 \\
\beta_1^3 & \beta_2^3 & 1 & 0 \\
\beta_1^4 & \beta_2^4 & \beta_3^4 & 1
\end{pmatrix}$. In this case, we denote $M^4(P, \lambda)$ by $M^4(\beta_1^2, \beta_1^3, \beta_2^3, \beta_1^4, \beta_2^4, \beta_3^4)$.

\begin{theorem} \label{bott_mfd_48}
 Let $\beta_1^2= 1$. If at least one of $\{ \beta_1^3, \beta_2^3\}$ is $1$, and at least two of $\{ \beta_1^4, \beta_2^4, \beta_3^4\}$ are $1$, then $8 \leq \mathbf{TC}(M^4(1, \beta_1^3, \beta_2^3, \beta_1^4, \beta_2^4, \beta_3^4)) \leq 9$.
\end{theorem}

\begin{proof}
The generators of the ideal $\Tilde{I}$ in the cohomology ring $H^*(M^4(1, \beta_1^3, \beta_2^3, \beta_1^4, \beta_2^4, \beta_3^4); \ZZ_2)$ are $\alpha_j= x_jy_j$ where $j=1, 2, 3, 4$. Let $\mathfrak{a}_j:= 1 \otimes y_j- y_j \otimes 1$ for $j= 1, 2, 3, 4$. Then $\mathfrak{a_j}$ is in the ideal of the zero-divisors of $H^*(M^4(1, \beta_1^3, \beta_2^3, \beta_1^4, \beta_2^4, \beta_3^4); \ZZ_2)$. Note that $\mathbf{TC}(M^4(1, \beta_1^3, \beta_2^3, \beta_1^4, \beta_2^4, \beta_3^4)) \leq 9$ by Proposition \ref{tc_small_cover_upper_lower}. By \cite[Theorem 5]{Naz11}, it is enough to consider the following manifolds to prove the claim; $M^4(1, 1, 0, 1, 1, 0), M^4(1, 0, 1, 1, 1, 0), M^4(1, 0, 1, 0, 1, 1), M^4(1, 0, 1, 1, 0, 1)$, and $M^4(1, 1, 1, 1, 1, 0)$
\begin{enumerate}
\item \label{bott_mfd_48_case1}
Consider the real Bott manifold $M^4(1, 1, 0, 1, 1, 0)$. Then from Proposition \ref{cohomoring_smallcover_polyt}, the ideal $\Tilde{J}$ is generated by the elements $x_1+ y_1, x_2+ y_1+ y_2, x_3+ y_1+ y_3$, and $x_4+ y_1+ y_2+ y_4$. So $x_1= y_1, x_2= y_1+ y_2, x_3= y_1+ y_3$, and $x_4= y_1+ y_2+ y_4$ in $H^*(M^4(1,1,0,1,1,0); \ZZ_2)$. Therefore, we have $y_1^2=0, y_2^2= y_1y_2, y_3^2= y_1y_3, y_4^2= y_1y_4+ y_2y_4$.
Now,
\begin{align*}
\mathfrak{a}_2\mathfrak{a}_3^3 \mathfrak{a}_4^3 &=  (1 \otimes y_2- y_2 \otimes 1) (1 \otimes y_3^3- y_3 \otimes y_3^2+ y_3^2 \otimes y_3- y_3^3 \otimes 1) (1 \otimes y_4^3- y_4 \otimes y_4^2+ y_4^2 \otimes y_4- y_4^3 \otimes 1) \nonumber \\
&=  y_1y_2y_3y_4 \otimes y_1y_2y_3- y_1y_2y_3 \otimes y_1y_2y_3y_4+ y_1y_2y_3y_4 \otimes y_1y_3y_4-  y_1y_3y_4 \otimes y_1y_2y_3y_4.\nonumber
\end{align*}
So, the product $\mathfrak{a}_2\mathfrak{a}_3^3 \mathfrak{a}_4^3$ contains an element $y_1y_2y_3 \otimes y_1y_2y_3y_4$ which is non-zero. Therefore the zero-divisors-cup-length of $H^*(M^4(1,1,0,1,1,0); \ZZ_2)$ is greater than or equal to $7$. Hence by Proposition \ref{zcl_less_tc_less_dim}, we have $8 \leq \mathbf{TC}(M^4(1,1,0,1,1,0))$.

\item
Consider the real Bott manifold $M^4(1,0,1,1,1,0)$. Then from Proposition \ref{cohomoring_smallcover_polyt}, the ideal $\Tilde{J}$ is generated by the elements $x_1+ y_1, x_2+ y_1+ y_2, x_3+ y_2+ y_3$, and $x_4+ y_1+ y_2+ y_4$. So, $x_1= y_1, x_2= y_1+ y_2, x_3= y_2+ y_3$, and $x_4= y_1+ y_2+ y_4$ in $H^*(M^4(1,0,1,1,1,0); \ZZ_2)$. Now,
\begin{align*}
\mathfrak{a}_2\mathfrak{a}_3^3 \mathfrak{a}_4^3 &=  (1 \otimes y_2- y_2 \otimes 1) (1 \otimes y_3^3- y_3 \otimes y_3^2+ y_3^2 \otimes y_3- y_3^3 \otimes 1) (1 \otimes y_4^3- y_4 \otimes y_4^2 \\
& + y_4^2 \otimes y_4- y_4^3 \otimes 1) \nonumber \\
&=  y_1y_2y_3y_4 \otimes (y_1y_2y_3+ y_2y_3y_4+ y_1y_3y_4)- (y_1y_2y_3+ y_2y_3y_4+ y_1y_3y_4) \otimes y_1y_2y_3y_4.\nonumber
\end{align*}
So, the product $\mathfrak{a}_2\mathfrak{a}_3^3 \mathfrak{a}_4^3$ contains an element $y_1y_2y_3 \otimes y_1y_2y_3y_4$ which is non-zero. Therefore the zero-divisors-cup-length of $H^*(M^4(1,0,1,1,1,0); \ZZ_2)$ is greater than or equal to $7$. Hence by Proposition \ref{zcl_less_tc_less_dim}, we have $8 \leq \mathbf{TC}(M^4(1,0,1,1,1,0))$.

\item
Consider the real Bott manifold $M^4(1,0,1,0,1,1)$. Then from Proposition \ref{cohomoring_smallcover_polyt}, the ideal $\Tilde{J}$ is generated by the elements $x_1+ y_1, x_2+ y_1+ y_2, x_3+ y_2+ y_3$, and $x_4+ y_2+ y_3+ y_4$. So, $x_1= y_1, x_2= y_1+ y_2, x_3= y_2+ y_3$, and $x_4= y_2+ y_3+ y_4$ in $H^*(M^4(1,0,1,0,1,1); \ZZ_2)$. Now,
\begin{align*}
\mathfrak{a}_2\mathfrak{a}_3^3 \mathfrak{a}_4^3 &=  (1 \otimes y_2- y_2 \otimes 1) (1 \otimes y_3^3- y_3 \otimes y_3^2+ y_3^2 \otimes y_3- y_3^3 \otimes 1) (1 \otimes y_4^3- y_4 \otimes y_4^2 \\
& + y_4^2 \otimes y_4- y_4^3 \otimes 1) \nonumber \\
&=  y_1y_2y_3y_4 \otimes (y_1y_2y_3+ y_2y_3y_4+ y_1y_2y_4)- (y_1y_2y_3+ y_2y_3y_4+ y_1y_2y_4) \otimes y_1y_2y_3y_4.\nonumber
\end{align*}
So, the product $\mathfrak{a}_2\mathfrak{a}_3^3 \mathfrak{a}_4^3$ contains an element $y_1y_2y_3 \otimes y_1y_2y_3y_4$ which is non-zero. Therefore the zero-divisors-cup-length of $H^*(M^4(1,0,1,0,1,1); \ZZ_2)$ is greater than or equal to $7$. Hence by Proposition \ref{zcl_less_tc_less_dim}, we have $8 \leq \mathbf{TC}(M^4(1,0,1,0,1,1))$.

\item
Consider the real Bott manifold $M^4(1,0,1,1,0,1)$. Then from Proposition \ref{cohomoring_smallcover_polyt}, the ideal $\Tilde{J}$ is generated by the elements $x_1+ y_1, x_2+ y_1+ y_2, x_3+ y_2+ y_3$, and $x_4+ y_1+ y_3+ y_4$. So $x_1= y_1, x_2= y_1+ y_2, x_3= y_2+ y_3$, and $x_4= y_1+ y_3+ y_4$ in $H^*(M^4(1,0,1,1,0,1); \ZZ_2)$. Now,
\begin{align*}
\mathfrak{a}_2\mathfrak{a}_3^3 \mathfrak{a}_4^3 &=  (1 \otimes y_2- y_2 \otimes 1) (1 \otimes y_3^3- y_3 \otimes y_3^2+ y_3^2 \otimes y_3- y_3^3 \otimes 1) (1 \otimes y_4^3- y_4 \otimes y_4^2 \\
& + y_4^2 \otimes y_4- y_4^3 \otimes 1) \nonumber \\
&=  y_1y_2y_3y_4 \otimes y_1y_2y_4- y_1y_2y_4 \otimes y_1y_2y_3y_4+ y_1y_2y_3y_4 \otimes y_1y_3y_4- y_1y_3y_4 \otimes y_1y_2y_3y_4.\nonumber
\end{align*}
So, the product $\mathfrak{a}_2\mathfrak{a}_3^3 \mathfrak{a}_4^3$ contains an element $y_1y_2y_4 \otimes y_1y_2y_3y_4$ which is non-zero. Therefore the zero-divisors-cup-length of $H^*(M^4(1,0,1,1,0,1); \ZZ_2)$ is greater than or equal to $7$. Hence by Proposition \ref{zcl_less_tc_less_dim}, we have $8 \leq \mathbf{TC}(M^4(1,0,1,1,0,1))$.

\item
Consider the real Bott manifold $M^4(1,1,1,1,1,0)$. Then from Proposition \ref{cohomoring_smallcover_polyt}, the ideal $\Tilde{J}$ is generated by the elements $x_1+ y_1, x_2+ y_1+ y_2, x_3+ y_1+ y_2+ y_3$, and $x_4+ y_1+ y_2+ y_4$. So, $x_1= y_1, x_2= y_1+ y_2, x_3= y_1+ y_2+ y_3$, and $x_4= y_1+ y_2+ y_4$ in $H^*(M^4(1,1,1,1,1,0); \ZZ_2)$. Now,
\begin{align*}
\mathfrak{a}_2\mathfrak{a}_3^3 \mathfrak{a}_4^3 &=  (1 \otimes y_2- y_2 \otimes 1) (1 \otimes y_3^3- y_3 \otimes y_3^2+ y_3^2 \otimes y_3- y_3^3 \otimes 1) (1 \otimes y_4^3- y_4 \otimes y_4^2 \\
& + y_4^2 \otimes y_4- y_4^3 \otimes 1) \nonumber \\
&=  y_1y_2y_3y_4 \otimes y_2y_3y_4 - y_2y_3y_4 \otimes y_1y_2y_3y_4.\nonumber
\end{align*}
So, the product $\mathfrak{a}_2 \mathfrak{a}_3^3 \mathfrak{a}_4^3$ contains an element $y_2y_3y_4 \otimes y_1y_2y_3y_4$ which is non-zero. Therefore the zero-divisors-cup-length of $H^*(M^4(1,1,1,1,1,0); \ZZ_2)$ is greater than or equal to $7$. Hence by Proposition \ref{zcl_less_tc_less_dim}, we have $8 \leq \mathbf{TC}(M^4(1,1,1,1,1,0))$.

\end{enumerate}
\end{proof}

\section{Symmetric topological complexity of small covers} \label{sec_symm_tc_small_cover}

In this section, we recall the definition of symmetric topological complexity. Then we compute this invariant for a class of small covers. 

Let $Y$ be a path-connected space. The path fibration ${ \pi \colon PY \rightarrow Y \times Y }$ restricts to a fibration
\begin{equation} \label{symm_second_eq}
 \pi' \colon P'Y \rightarrow F(Y;2),   
\end{equation}
 
\noindent where $F(Y;2)=\lbrace (x,y) \in Y \times Y ~|~ x \neq y \rbrace$ is the space of ordered pairs of distinct points in $Y$, and $P'Y$ is the subspace $\lbrace \gamma \colon I \rightarrow Y ~|~ \gamma(0) \neq \gamma(1) \rbrace \subseteq PY$ consisting of paths with distinct endpoints.

 The group $\mathbb{Z}_2$ acts on $F(Y;2)$ by permutation of factors, and acts on $P'Y$ by sending a path $\gamma$ to its inverse $\bar{\gamma}$ given by $\bar{\gamma}(t)= \gamma (1-t)$. So, the group $\mathbb{Z}_2$ acting on the spaces $P'Y$ and $F(Y;2)$ freely. Observe that
${\pi' \colon P'Y \rightarrow F(Y;2)}$ 
is an equivariant map of free $\mathbb{Z}_2$-spaces. So, it induces a map
\begin{equation} \label{symm_third_eq}
  \pi'' \colon P'Y/\mathbb{Z}_2 \rightarrow B(Y;2), 
\end{equation}
where $B(Y;2)$ denotes the orbit space $F(Y;2)/\mathbb{Z}_2$ of unordered pairs of distinct points in $Y$. This map is also a fibration.
\begin{definition}
The symmetric topological complexity of $Y$, denoted by $\mathbf{TC}^S(Y)$, is defined to be one plus the sectional category of the fibration $\pi''$. In other words, $\textbf{TC}^S(Y)=1+ \mbox{secat}(\pi'')$.
\end{definition}
We adopt the convention that the sectional category of $p \colon E \rightarrow B$ vanishes if and only if $E= B= \emptyset$. The space $B(Y;2)$ is empty if and only if $Y$ is a single point, and so in this case, $\mathbf{TC}^S(Y)=1$. If $Y$ contains more than one point then $\mbox{secat}(\pi'') \geq 1$, and therefore $\mathbf{TC}^S(Y) \geq 2$.



\begin{example}
 Let $Y$ is a contractible space. Then there exists a continuous map $y \mapsto \gamma_y \in PY$ such that  $\gamma_y(0)= y$ and $\gamma_y(1)= y_0$. Then setting $s(a, b)$ to be equal to the concatenation of $\gamma_a$ and the inverse path to $\gamma_b$ gives a symmetric equivariant section of \eqref{symm_second_eq}. Therefore for any contractible space $Y$ with more than one point, we have $\mathbf{TC}^S(Y)= 2$. We note that if $Y$ is a path-connected space with $\mathbf{TC}^S(Y)= 2$, then $Y$ is contractible.
\end{example}

Let $N_Y$ be the sub-ring of $H^*(Y) \otimes H^*(Y)$ spanned by the norm elements (i.e., the elements of the form $x \otimes y+ y \otimes x$ with $x \neq y$). The following result follows from Corollary $9$, Proposition $10$, and Theorem $17$ in \cite{FG07}.

\begin{proposition} \label{lower_upper_bound_symm_tc}
    Let $Y$ be a closed smooth manifold. Then 
    $$ \mbox{max} \{ \mathbf{TC}(Y), \mbox{cl}(N_Y)+ 2 \} \leq \mathbf{TC}^S(Y) \leq 2 \mbox{dim} Y+ 1. $$
\end{proposition}

Next, we calculate the symmetric topological complexity of the circle.




\begin{corollary} \label{prop_rp1}
 If $P$ is an $1$-simplex, then $M^1(P, \lambda)= \mathbb{RP}^1$ and $\mathbf{TC}^S(\mathbb{RP}^1)= 3$.  
\end{corollary}
\begin{proof}
 Note that the non-zero element $1 \otimes y_1+ y_1 \otimes 1$ is the norm element of $H^*(M^1(P, \lambda); \ZZ_2) \otimes H^*(M^1(P, \lambda); \ZZ_2)$. So, by Proposition \ref{lower_upper_bound_symm_tc}, we get $3 \leq \mathbf{TC}^S(M^1(P, \lambda))$. The small cover over a $1$-simplex is $\mathbb{RP}^1= \mathbb{S}^1$. Therefore, by Proposition \ref{lower_upper_bound_symm_tc}, we get $$\mathbf{TC}^S(M^1(P, \lambda))\leq 2\mbox{dim}(M^1(P, \lambda))+ 1= 3.$$    
\end{proof}

We note that the conclusion of Corollary \ref{prop_rp1} can be obtained from \cite[Corollary 18]{FG07}. 

\begin{remark} \label{cup_length_zcl}
The element $1 \otimes y_j+ y_j \otimes 1$ is same as $1 \otimes y_j- y_j \otimes 1$ in $N_{M^n(P, \lambda)}$. Thus the zero-divisors-cup-length of $M^n(P, \lambda)$ is the same as the cup-length of $N_{M^n(P, \lambda)}$.  
\end{remark}

\begin{theorem}\label{symmtc_small_cover_special}
Let $M^n(P, \lambda)$ be a small cover other than $\mathbb{RP}^{n_1} \times \mathbb{RP}^{n_2}$ over $P= \Delta^{n_1} \times \Delta^{n_2}$.
\begin{enumerate}
 \item \label{SymTC_main_proof_case1}
 Let $n_2 \in \mathcal{S}$ with $n_2 > n_1$. Then $2^{r_1}+2^{r_2} \leq \mathbf{TC}^S(M^n(P, \lambda))$.
 \item \label{SymTC_main_proof_case2}
 Let $n_2 \in \mathcal{S}$ with $n_2$ divides $n_1$. Then $2^r+1 \leq \mathbf{TC}^S(M^n(P, \lambda))$.
\item \label{SymTC_main_proof_case3}
Let $n_2 \in \mathcal{S}+1$ with $n_2 > n_1+1$. Then, $ 2^r+1 \leq \mathbf{TC}^S(M^n(P, \lambda))$.

\item \label{SymTC_main_proof_case4}
Let $n_2 \in \mathcal{S}+2$ with $n_2 > n_1+2$. Then, $2^r+1 \leq \mathbf{TC}^S(M^n(P, \lambda))$.
\end{enumerate}
In particular, if $n=2^{s-1}$, then for the cases \eqref{SymTC_main_proof_case2}, \eqref{SymTC_main_proof_case3}, and \eqref{SymTC_main_proof_case4}, we have $\mathbf{TC}^S(M^n(P, \lambda))= 2n+1$.
\end{theorem}
\begin{proof}
 Let $N_{M^n(P, \lambda)}$ denotes the sub-ring of $H^*(M^n(P, \lambda); \ZZ_2) \otimes H^*(M^n(P, \lambda); \ZZ_2)$ spanned by the norm elements. Consider the norm elements $\mathfrak{a}_j: =1 \otimes y_j+ y_j \otimes 1$ in  $N_{M^n(P, \lambda)}$ for $j=1, 2$. Then the proof follows from Theorem \ref{tc_small_cover_special}, Proposition \ref{lower_upper_bound_symm_tc}, and Remark \ref{cup_length_zcl}.







\end{proof}

\begin{theorem}
For $n \geq 3$, let the elements $\beta_k^{k+1}$ in the Bott matrix \eqref{lower_triangular_matrix} be $1$ for $k=1,..., n-1$, and the remaining elements $\beta_l^m$ be zero for $l= 1,...,n-2$ and $m= 3,...,n$. If $n \leq 2^r-1 < 2n$, then the symmetric topological complexity of the real Bott manifold $M^n(P, \lambda)$ is greater than or equal to $2^r+1$. In particular, if $n= 2^{s-1}$, then
$\mathbf{TC}^S(M^n(P, \lambda))= 2n+1$.
\end{theorem}
\begin{proof}
The proof follows from Theorem \ref{real_bott_tc_matrix_change}, Proposition \ref{lower_upper_bound_symm_tc}, and Remark \ref{cup_length_zcl}.  
\end{proof}



\begin{remark} \hfill
\begin{enumerate}
\item From Theorem \ref{bott_mfd_36}, Proposition \ref{lower_upper_bound_symm_tc}, and Remark \ref{cup_length_zcl}, we get, 
$${\mathbf{TC}^S(M^3(1, 1, 0))= \mathbf{TC}^S(M^3(1, 0, 1))= \mathbf{TC}^S(M^3(1, 1, 1))= 7}.$$

 \item Let $\beta_1^2= 1$. If at least one of $\{ \beta_1^3, \beta_2^3\}$ is $1$, and at least two of $\{ \beta_1^4, \beta_2^4, \beta_3^4\}$ are $1$. Then from Theorem \ref{bott_mfd_48}, Proposition \ref{lower_upper_bound_symm_tc}, and Remark \ref{cup_length_zcl}, we get, $${\mathbf{TC}^S(M^4(1, \beta_1^3, \beta_2^3, \beta_1^4, \beta_2^4, \beta_3^4))= 9}.$$
\end{enumerate}
\end{remark}

\section{$\mathcal{D}$-topological complexity of small covers} \label{sec_dtc_small_cover}

In this section, we recall the $\mathcal{D}$-topological complexity and the LS one-category of a space. We compute LS one-category for all real Bott manifolds and for a class of small covers. Then, we give some bounds for the $\mathcal{D}$-topological complexity of small covers.

\begin{definition}
Let $Y$ be a path-connected space with the fundamental group $G= \pi_1(Y,y_0)$. The $\mathcal{D}$-topological complexity, denoted by $\mathbf{TC}^{\mathcal{D}}(Y)$, is defined as the minimal number $k$ such that $Y \times Y$ can be covered by $k$ open subsets $U_1,..., U_k$ with the property that for each $i \in \{1,...,k \}$ and for every choice of the base point $u_i \in U_i$, the homomorphism $\pi_1(U_i,u_i) \rightarrow \pi_1(Y \times Y,u_i)$ induced by the inclusion $U_i \rightarrow Y \times Y$ takes values in a subgroup conjugate to the diagonal $\Delta \subseteq G \times  G$.
\end{definition}
Note that there is an isomorphism $\pi_1(Y \times Y,u_i) \rightarrow \pi_1(Y \times Y,(y_0, y_0)) \cong G \times G$ determined uniquely up to conjugation, and the diagonal inclusion $Y \rightarrow Y \times Y$ induces the inclusion $G \rightarrow G \times G$ onto the diagonal $\Delta$.





We recall the Lusternik-Schnirelmann one-category (in short LS one-category) of a space which is denoted by $\mbox{cat}_1(Y)$ for a space $Y$.

\begin{definition}
Let $Y$ be a connected, locally path-connected, and semi-locally simply connected space with the universal cover $p \colon \widetilde{Y} \rightarrow Y$. Then the LS one-category is the sectional category of the map $p$. That is, $\mbox{cat}_1(Y)= \mbox{secat}(p)$. 
\end{definition}


Similar to $\mbox{cat}(Y)$ and $\mathbf{TC}(Y)$ there is a relation between $\mbox{cat}_1(Y)$ and $\mathbf{TC}^{\mathcal{D}}(Y)$.

\begin{proposition} \cite[Proposition 2.4, Proposition 2.11]{FGLO19} \label{cat1_less_TCD_less_tc}
If $Y$ is a connected, locally path-connected, and semi-locally simply connected topological space, then
$$\mbox{cat}_1(Y) \leq \mathbf{TC}^{\mathcal{D}}(Y) \leq \mbox{min} \{ \mathbf{TC}(Y), \mbox{cat}_1(Y \times Y)\}.$$
\end{proposition}

We recall a result that gives a lower bound for the sectional category of fibrations.

\begin{proposition} \cite[Proposition 9.14]{CLOT03} \label{secat_one_cat}
   Let $F \rightarrow E \xrightarrow{p} B$ be a fibration. If there exists $y_1,..., y_k \in H^*(B;R)$ with $p^*(y_1)= \cdots= p^*(y_k)= 0$ and $y_1 \cup \cdots \cup y_k \neq 0$, then $\mbox{secat}(p) \geq k+1$.
   \end{proposition}




The following result gives the computation of LS one-category of infinitely many small covers.

\begin{theorem} \label{one_cat_small_cover_moment_mfold}
 Let $M^n(P, \lambda)$ be a small cover over a simple polytope $P$ such that $\mathbb{R}\mathcal{Z}_{K_P}$ is simply connected. Then $\mbox{cat}_1(M^n(P, \lambda))= n+1$.
\end{theorem}

\begin{proof}
    
Consider the principal $\mathbb{Z}_2^m$-bundle map $p \colon \mathbb{R}\mathcal{Z}_{K_P} \rightarrow M^n(P; \lambda)$ given by Proposition \ref{small_cover_moment_manifold}. So we get the induced graded ring homomorphism
$$p^* \colon H^*(M^n(P, \lambda); \mathbb{Z}_2) \rightarrow H^*(\mathbb{R}\mathcal{Z}_{K_P}; \mathbb{Z}_2).$$
 
 Note that $p^*$ carries $H^j(M^n(P, \lambda); \mathbb{Z}_2)$ to $H^j(\mathbb{R}\mathcal{Z}_{K_P}; \mathbb{Z}_2)$. Now, for $j=1$,
 $H^1(\mathbb{R}\mathcal{Z}_{K_P}; \mathbb{Z}_2)=0$ as $\mathbb{R}\mathcal{Z}_{K_P}$ is simply connected. Therefore, each $v$ in $H^1(M^n(P, \lambda); \mathbb{Z}_2)$ maps to $0$ in $H^*(\mathbb{R}\mathcal{Z}_{K_P}; \mathbb{Z}_2)$. Hence $p^*(v)= 0$ for $v \in H^1(M^n(P, \lambda); \mathbb{Z}_2)$. Since $P$ is a simple polytope, at each vertex, exactly $n$ many facets intersect. So the cup product of corresponding $n$ indeterminates is non-zero. Therefore, by Proposition \ref{secat_one_cat}, $\mbox{secat}(p) \geq n+1$. Since $p$ is the universal cover, so by definition of LS one-category, $\mbox{secat}(p)= \mbox{cat}_1(M^n(P, \lambda))$. Therefore $n+ 1 \leq \mbox{cat}_1(M^n(P, \lambda))$.

On the other hand, we know that $\mbox{cat}_1(M^n(P, \lambda)) \leq \mbox{cat}(M^n(P, \lambda))$ as discussed in \cite{FGLO19}. Since $\mbox{cat}(M^n(P, \lambda))= n+1$ (by Theorem \ref{small_cover_lscat}), so $\mbox{cat}_1(M^n(P, \lambda)) \leq n+1$. Hence, we get the result.
\end{proof}

\begin{corollary} \label{one_cat_nj_2}
 Let $M^n(P, \lambda)$ be a small cover over $P= \prod_{j=1}^m \Delta^{n_j}$ such that $n_j \geq 2$ for $j= 1,..., m$. Then $\mbox{cat}_1(M^n(P, \lambda))= n+1$.
\end{corollary}

\begin{proof}
The moment angle manifold $\mathbb{R}\mathcal{Z}_{K_P}$ for the polytope $P= \prod_{j=1}^m \Delta^{n_j}$ is $\mathbb{S}^{n_1} \times \cdots \times \mathbb{S}^{n_m}$. Thus, $\mathbb{R}\mathcal{Z}_{K_P}$ is simply connected and orientable for $n_j \geq 2$ for $j= 1,..., m$. Therefore, by Theorem \ref{one_cat_small_cover_moment_mfold}, $\mbox{cat}_1(M^n(P, \lambda))= n+ 1$.
\end{proof}

\begin{theorem} \label{one_cat_small_cover_exact}
 Let $M^n(P, \lambda)$ be a small cover over $P= \prod_{j=1}^m \Delta^{n_j}$. Then $\mbox{cat}_1(M^n(P, \lambda))= n+1$.
\end{theorem}

\begin{proof}
Corollary \ref{one_cat_nj_2} gives the proof for all $n_j \geq 2$.

Now consider the small cover $M^n(P, \lambda)$ over $P= \prod_{j=1}^m \Delta^{n_j}$ where some $n_j=1$. Without loss of generality, we assume that $n_1=n_2= \cdots=n_{s-1}=1$ and the remaining $n_j$'s are greater than or equal to $2$. Then the map 
$\Bar{p} \colon (\prod_1^{s-1} \mathbb{R} \times \prod_{j=s}^m \mathbb{S}^{n_j}) \rightarrow  (\prod_1^{s-1}\mathbb{S}^1 \times \prod_{j=s}^m \mathbb{S}^{n_j})/ \mathbb{Z}_2^m= M^n(P, \lambda)$
is the universal cover where $\RR \rightarrow \mathbb{S}^1$ is given by exponential map. This induces a ring homomorphism 
$$\Bar{p}^* \colon H^*(M^n(P, \lambda); \mathbb{Z}_2) \rightarrow H^*(\prod_1^{s-1} \mathbb{R} \times \prod_{j=s}^m \mathbb{S}^{n_j}; \mathbb{Z}_2).$$

We know that the cohomology ring of $M^n(P, \lambda)$ is generated by $y_1,..., y_m$ and $y_j^{n_j} \neq 0$ for $j=1,2,..., m$. Hence 
\begin{equation*} 
    p^*(y_1)= \cdots= p^*(y_m)= 0.
 \end{equation*}

 Since all $y_j^{n_j} \neq 0$ for $j= 1,..., m$, so
 \[
  y_1 \cup \cdots \cup y_{s-1} \cup \underbrace{y_s \cup \cdots \cup y_s }_\text{$n_s$ times} \cup \cdots \cup
    \underbrace{y_m \cup \cdots \cup y_m}_\text{$n_m$ times} \neq 0.
\]
  Therefore, by Proposition \ref{secat_one_cat}, $\mbox{secat}(\Bar{p}) \geq n_1+ \cdots+ n_m+ 1=n+1$. By the definition of LS one-category, $\mbox{secat}(p)= \mbox{cat}_1(M^n(P, \lambda))$. Therefore $n+1 \leq \mbox{cat}_1(M^n(P, \lambda))$. Using Theorem \ref{small_cover_lscat}, $\mbox{cat}(M^n(P, \lambda))= n+1$. Thus $\mbox{cat}_1(M^n(P, \lambda)) \leq n+1$. Hence $\mbox{cat}_1(M^n(P, \lambda))=n+1$.
\end{proof}

We give some bounds on $\mathbf{TC}^{\mathcal{D}}(M^n(P, \lambda))$ in the following.

\begin{theorem}
 Let $M^n(P, \lambda)$ be a small cover over a product of simplices $P$. Then $$ n+1 \leq \mathbf{TC}^{\mathcal{D}}(M^n(P, \lambda)) \leq 2n+1.$$ In particular, if $M^n(P, \lambda)= \mathbb{RP}^{n_1} \times \cdots \times \mathbb{RP}^{n_m}$ with $n_j \in \{ 1, 3, 7 \}$, then $\mathbf{TC}^{\mathcal{D}}(M^n(P, \lambda))= n+1$.
\end{theorem}

\begin{proof}
By Proposition \ref{cat1_less_TCD_less_tc}, we have $\mbox{cat}_1(M^n(P, \lambda)) \leq \mathbf{TC}^{\mathcal{D}}(M^n(P, \lambda))$. Therefore, by Theorem \ref{one_cat_small_cover_exact}, we have $n+1 \leq \mathbf{TC}^{\mathcal{D}}(M^n(P, \lambda))$.
By Proposition \ref{cat1_less_TCD_less_tc}, $\mathbf{TC}^{\mathcal{D}}(M^n(P, \lambda)) \leq \mathbf{TC}(M^n(P, \lambda))$. So, the upper bound of $\mathbf{TC}^{\mathcal{D}}(M^n(P, \lambda))$ is $2n+1$, i.e., $\mathbf{TC}^{\mathcal{D}}(M^n(P, \lambda)) \leq 2n+ 1$.

The second part follows from Corollary \ref{tc_small_cover_special_case} and $\mathbf{TC}^{\mathcal{D}}(M^n(P, \lambda)) \leq \mathbf{TC}(M^n(P, \lambda))$. In this case $\mathbf{TC}^{\mathcal{D}}(M^n(P, \lambda)) \leq n+1$.
\end{proof}

{\bf Acknowledgment}. The authors thank Navnath for the helpful discussion. The first author would like to thank the Indian Institute of Technology Madras for the Ph.D. fellowship. The second author would like to thank UGC India for his Ph.D. fellowship. The last author would like a thank ISI Kolkata for its financial support.

\bibliographystyle{amsalpha}
\bibliography{bib}

\end{document}